\documentclass[journal]{IEEEtran}

\usepackage[utf8]{inputenc}
\usepackage{color}
\usepackage{array}
\usepackage{verbatim}
\usepackage{float}
\usepackage{amsmath}
\usepackage{amsthm}
\usepackage{amssymb}
\usepackage{graphicx}
\usepackage{longtable}
\usepackage[unicode=true,
 bookmarks=false,
 breaklinks=false,pdfborder={0 0 1},colorlinks=false]
 {hyperref}
\hypersetup{
 colorlinks,bookmarksopen,bookmarksnumbered,citecolor=blue,urlcolor=blue}
\usepackage{cite}
\makeatletter
 \let\oldforeign@language\foreign@language
 \DeclareRobustCommand{\foreign@language}[1]{%
   \lowercase{\oldforeign@language{#1}}}

 \let\oldforeign@language\foreign@language
 \DeclareRobustCommand{\foreign@language}[1]{%
   \lowercase{\oldforeign@language{#1}}}

\ifCLASSINFOpdf
\else
\fi

\newtheorem{defn}{Definition}

\newtheorem{lem}{Lemma}

\newtheorem{thm}{Theorem}
\newtheorem{rem}{Remark}

\begin{document}
\bstctlcite{IEEEexample:BSTcontrol}
\onecolumn
\noindent\rule{18.1cm}{2pt}\\
\underline{To cite this article:}
{\bf{\textcolor{red}{H. A. Hashim, L. J. Brown, and K. McIsaac, "Nonlinear Explicit Stochastic Attitude Filter on SO(3),"  in Proceedings of the 57th IEEE Conference on Decision and Control (CDC), 2018, pp. 1210 - 1216.}}}\\
\noindent\rule{18.1cm}{2pt}\\


\vspace{40pt}\noindent Please note that where the full-text provided is the Author Accepted Manuscript or Post-Print version this may differ from the final Published version. { \bf To cite this publication, please use the final published version.}\\

\textbf{
	\begin{center}
		Personal use of this material is permitted. Permission from the author(s) and/or copyright holder(s), must be obtained for all other uses, in any current or future media, including reprinting or republishing this material for advertising or promotional purposes.\vspace{60pt}\\
	\end{center}
\vspace{400pt}
}
\footnotesize{ \bf
	\vspace{20pt}\noindent Please contact us and provide details if you believe this document breaches copyrights. We will remove access to the work immediately and investigate your claim.
} 

\normalsize

\twocolumn
\title{Nonlinear Explicit Stochastic Attitude Filter on SO(3)}

\author{Hashim~A.~Hashim, Lyndon J. Brown, and~Kenneth McIsaac
\thanks{H. A. Hashim, L. J. Brown and K. McIsaac are with the Department of Electrical and Computer Engineering,
University of Western Ontario, London, ON, Canada, N6A-5B9, e-mail: hmoham33@uwo.ca, lbrown@uwo.ca and kmcisaac@uwo.ca.}
}


\markboth{--,~Vol.~-, No.~-, \today}{Hashim \MakeLowercase{\textit{et al.}}: Nonlinear Stochastic Attitude Filter on the Special Orthogonal Group}
\markboth{}{Hashim \MakeLowercase{\textit{et al.}}: Nonlinear Stochastic Attitude Filter on the Special Orthogonal Group}

\maketitle

\begin{abstract}
This work proposes a nonlinear stochastic filter evolved on the Special
Orthogonal Group $\mathbb{SO}\left(3\right)$ as a solution to the
attitude filtering problem. One of the most common potential functions
for nonlinear deterministic attitude observers is studied and reformulated
to address the noise attached to the attitude dynamics. The resultant
estimator and correction factor demonstrate convergence properties
and remarkable ability to attenuate the noise. The stochastic dynamics
of the attitude problem are mapped from $\mathbb{SO}\left(3\right)$
to Rodriguez vector. The proposed stochastic filter evolved on $\mathbb{SO}\left(3\right)$
guarantees that errors in the Rodriguez vector and estimates steer
very close to the neighborhood of the origin and that the errors are
semi-globally uniformly ultimately bounded in mean square. Simulation
results illustrate the robustness of the proposed filter in the presence
of high uncertainties in measurements.  
\end{abstract}

%

\IEEEpeerreviewmaketitle{}

\section{Introduction}

The orientation of a rigid-body is termed attitude, and attitude estimation
is an essential subtask in robotics applications \cite{mahony2008nonlinear,crassidis2007survey}.
Unfortunately, the attitude cannot be accurately measured, however,
the available measurements from sensors attached to the body-frame
and inertial-frame coupled with an attitude filter allow reasonably
accurate estimation of the true attitude. The moving vehicles are
normally equipped with low-cost inertial measurement units (IMUs)
which are very sensitive to noise and bias components, complicating
the attitude estimation \cite{mahony2008nonlinear,crassidis2007survey}.

Historically, the attitude filtering problem has been addressed using
Gaussian filters based mainly on the structure of the Kalman filter
(KF) \cite{crassidis2007survey}. The family of attitude Gaussian
filters includes KF \cite{choukroun2006novel}, extended KF (EKF)
\cite{lefferts1982kalman}, multiplicative EKF (MEKF) \cite{markley2003attitude},
and others. However, attitude Gaussian filters have proven to be inefficient
if the vehicle is equipped with low quality sensors \cite{mahony2008nonlinear,crassidis2007survey}.
Other filtering techniques such as Unscented KF (UKF) \cite{crassidis2003unscentedy}
and particle filters (PFs) \cite{Arulampalam2002Particle} provide
a more precise estimation even when the low-quality sensors are used.
However, the computational cost of the above-mentioned filters is
higher \cite{crassidis2007survey}. It should be remarked that the
Gaussian filters in \cite{choukroun2006novel,lefferts1982kalman,markley2003attitude}
as well as UKF \cite{crassidis2003unscentedy}
and PFs \cite{Arulampalam2002Particle} are quaternion based which
does not provide a unique representation of the attitude \cite{shuster1981three}.

The deficiencies of Gaussian filters, UKF and PFs, in addition to
the development of low-cost IMUs, motivated researchers to design
nonlinear deterministic attitude filters, such as \cite{mahony2008nonlinear,zlotnik2017nonlinear,grip2012attitude}.
These filters have better tracking performance than Gaussian filters
\cite{mahony2008nonlinear} and require less computational power when
compared with UKF and PFs \cite{crassidis2007survey}. In addition,
nonlinear deterministic attitude filters evolve directly on $\mathbb{SO}\left(3\right)$
which is nonsingular in parameterization and has a unique representation.
The deterministic filters proposed in \cite{mahony2008nonlinear,zlotnik2017nonlinear,grip2012attitude}
can be easily fitted given two or more vectorial measurements and
a rate gyroscope measurement, however, the selected potential functions
in \cite{mahony2008nonlinear,zlotnik2017nonlinear,grip2012attitude}
were kept unchanged. The potential function in \cite{mahony2008nonlinear,zlotnik2017nonlinear,grip2012attitude}
fits nonlinear deterministic attitude filters on $\mathbb{SO}\left(3\right)$
assuming that the rate gyro measurements are corrupted only with constant
bias and are noise free. However, the environment is noisy \cite{hashim2017adaptive,hashim2017neuro}
and the kinematics of the nonlinear attitude problem on $\mathbb{SO}\left(3\right)$
in its natural stochastic sense need to be considered.

The main challenge is that the attitude problem is 1) modeled on the
Lie group of $\mathbb{SO}\left(3\right)$ which is nonlinear and 2)
the attitude dynamics are a function of angular velocity measurements
which are corrupted with noise components. Therefore, the randomness
and uncertain behavior in attitude kinematics prompted the proposal
of nonlinear stochastic attitude filter on $\mathbb{SO}\left(3\right)$
based on the selection of a new potential function. Hence, in the
case where angular velocity measurements are contaminated with noise,
the stochastic filter would be able to guarantee that, %
{} 1) the error is regulated to an arbitrarily small neighborhood of
the equilibrium point in probability; and 2) the error is semi-globally
uniformly ultimately bounded (SGUUB) in mean square.

The rest of the paper is organized as follows: Section \ref{sec:SO3STCH_EXP_Math-Notations}
gives an overview of mathematical notation and preliminaries. The
problem is formulated in stochastic sense in Section \ref{sec:SO3STCH_EXP_Problem-Formulation-in}.
The nonlinear stochastic filter on $\mathbb{SO}\left(3\right)$ is
proposed and the stability analysis is presented in Section \ref{sec:SO3STCH_EXP_Stochastic-Complementary-Filters}.
Section \ref{sec:SO3STCH_EXP_Simulation} demonstrates the numerical
results. Finally, closing notes are provided in Section \ref{sec:SO3STCH_EXP_Conclusion}.

\section{Math Notation \label{sec:SO3STCH_EXP_Math-Notations}}

In this paper, $\mathbb{R}^{n}$ is the real $n$-dimensional space
while $\mathbb{R}^{n\times m}$ denotes the real $n\times m$ dimensional
space. For $x\in\mathbb{R}^{n}$, the Euclidean norm is defined as
$\left\Vert x\right\Vert =\sqrt{x^{\top}x}$, where $^{\top}$ is
the transpose of a component. $\mathcal{C}^{n}$ denotes the set of
functions with continuous $n$th partial derivatives. $\mathbb{P}\left\{ \cdot\right\} $,
$\mathbb{E}\left[\cdot\right]$, ${\rm exp}\left(\cdot\right)$, and
${\rm Tr}\left\{ \cdot\right\} $ refer to probability, expected value,
exponential, and trace of a component, respectively. $\lambda\left(\cdot\right)$
is the set of eigenvalues of the associated matrix while $\underline{\lambda}\left(\cdot\right)$
is the minimum singular value. $\mathbf{I}_{n}$ denotes identity
with dimensions $n$-by-$n$, and $\underline{\mathbf{0}}_{n}\in\mathbb{R}^{n}$
is a zero column vector. $\mathbb{SO}\left(3\right)$ denotes the
Special Orthogonal Group, and the attitude of a rigid-body is defined
as a rotational matrix $R$:
\[
\mathbb{SO}\left(3\right):=\left\{ \left.R\in\mathbb{R}^{3\times3}\right|R^{\top}R=\mathbf{I}_{3}\text{, }{\rm det}\left(R\right)=1\right\} 
\]
where ${\rm det\left(\cdot\right)}$ is the determinant of the associated
matrix. The Lie-algebra of $\mathbb{SO}\left(3\right)$ is known as
$\mathfrak{so}\left(3\right)$ and is given by
\[
\mathfrak{so}\left(3\right):=\left\{ \left.\mathcal{Y}\in\mathbb{R}^{3\times3}\right|\mathcal{Y}^{\top}=-\mathcal{Y}\right\} 
\]
with $\mathcal{Y}$ being the space of skew-symmetric matrices. Define
the map $\left[\cdot\right]_{\times}:\mathbb{R}^{3}\rightarrow\mathfrak{so}\left(3\right)$
such that
\[
\mathcal{Y}=\left[y\right]_{\times}=\left[\begin{array}{ccc}
0 & -y_{3} & y_{2}\\
y_{3} & 0 & -y_{1}\\
-y_{2} & y_{1} & 0
\end{array}\right],\hspace{1em}y=\left[\begin{array}{c}
y_{1}\\
y_{2}\\
y_{3}
\end{array}\right]
\]
For all $\psi,\beta\in\mathbb{R}^{3}$, we have $\left[\psi\right]_{\times}\beta=\psi\times\beta$
where $\times$ is the cross product between the two vectors. Let
the vex operator be the inverse of $\left[\cdot\right]_{\times}$,
denoted by $\mathbf{vex}:\mathfrak{so}\left(3\right)\rightarrow\mathbb{R}^{3}$
such that $\mathbf{vex}\left(\mathcal{B}\right)=\beta$ for all $\beta\in\mathbb{R}^{3}$
and $\mathcal{B}\in\mathfrak{so}\left(3\right)$. Let $\boldsymbol{\mathcal{P}}_{a}$
denote the anti-symmetric projection operator on the Lie-algebra $\mathfrak{so}\left(3\right)$,
defined by $\boldsymbol{\mathcal{P}}_{a}:\mathbb{R}^{3\times3}\rightarrow\mathfrak{so}\left(3\right)$
such that $\boldsymbol{\mathcal{P}}_{a}\left(\mathcal{A}\right)=\frac{1}{2}\left(\mathcal{A}-\mathcal{A}^{\top}\right)\in\mathfrak{so}\left(3\right)$
for all $\mathcal{A}\in\mathbb{R}^{3\times3}$. The normalized Euclidean
distance of a rotation matrix on $\mathbb{SO}\left(3\right)$ is given
by the following equation
\begin{equation}
\left\Vert R\right\Vert _{I}:=\frac{1}{4}{\rm Tr}\left\{ \mathbf{I}_{3}-R\right\} \label{eq:SO3STCH_EXPL_Ecul_Dist}
\end{equation}
where $\left\Vert R\right\Vert _{I}\in\left[0,1\right]$. The attitude
of a rigid body can be constructed knowing angle of rotation $\alpha\in\mathbb{R}$
and axis parameterization $u\in\mathbb{R}^{3}$. The mapping of angle-axis
parameterization to $\mathbb{SO}\left(3\right)$ is defined by $\mathcal{R}_{\alpha}:\mathbb{R}\times\mathbb{R}^{3}\rightarrow\mathbb{SO}\left(3\right)$
such that
\begin{align}
\mathcal{R}_{\alpha}\left(\alpha,u\right) & =\mathbf{I}_{3}+\sin\left(\alpha\right)\left[u\right]_{\times}+\left(1-\cos\left(\alpha\right)\right)\left[u\right]_{\times}^{2}\label{eq:SO3STCH_EXPL_att_ang}
\end{align}
The attitude could be obtained knowing Rodriguez parameters vector
$\rho\in\mathbb{R}^{3}$ and $\mathcal{R}_{\rho}:\mathbb{R}^{3}\rightarrow\mathbb{SO}\left(3\right)$
\cite{shuster1993survey} is
\begin{align}
\mathcal{R}_{\rho}\left(\rho\right)= & \frac{1}{1+\left\Vert \rho\right\Vert ^{2}}\left(\left(1-\left\Vert \rho\right\Vert ^{2}\right)\mathbf{I}_{3}+2\rho\rho^{\top}+2\left[\rho\right]_{\times}\right)\label{eq:SO3STCH_EXPL_SO3_Rodr}
\end{align}
with direct substitution of \eqref{eq:SO3STCH_EXPL_SO3_Rodr} in \eqref{eq:SO3STCH_EXPL_Ecul_Dist}
for $\mathcal{R}_{\rho}=\mathcal{R}_{\rho}\left(\rho\right)$ one
obtains
\begin{equation}
||R||_{I}=\frac{1}{4}{\rm Tr}\left\{ \mathbf{I}_{3}-\mathcal{R}_{\rho}\right\} =\frac{||\rho||^{2}}{1+||\rho||^{2}}\label{eq:SO3STCH_EXPL_TR2}
\end{equation}
Likewise, the anti-symmetric projection operator of attitude $R$
in \eqref{eq:SO3STCH_EXPL_SO3_Rodr} can be defined as
\begin{align*}
\boldsymbol{\mathcal{P}}_{a}\left(R\right)=\frac{1}{2}\left(\mathcal{R}_{\rho}-\mathcal{R}_{\rho}^{\top}\right)= & 2\frac{1}{1+||\rho||^{2}}\left[\rho\right]_{\times}
\end{align*}
and the vex operator of the above-mentioned result is
\begin{equation}
\Phi\left(R\right)=\mathbf{vex}\left(\boldsymbol{\mathcal{P}}_{a}\left(R\right)\right)=2\frac{\rho}{1+||\rho||^{2}}\label{eq:SO3STCH_EXPL_VEX_Pa}
\end{equation}
where $\Phi\left(\cdot\right)$ is the composition mapping such that
$\Phi:=\mathbf{vex}\circ\boldsymbol{\mathcal{P}}_{a}$. The following
identities will be used in the subsequent derivations
\begin{align}
\left[\psi\times\beta\right]_{\times} & =\beta\psi^{\top}-\psi\beta^{\top},\quad\psi,\beta\in{\rm \mathbb{R}}^{3}\label{eq:SO3STCH_EXPL_Identity1}\\
\left[R\beta\right]_{\times} & =R\left[\beta\right]_{\times}R^{\top},\quad R\in\mathbb{SO}\left(3\right),\beta\in\mathbb{R}^{3}\label{eq:SO3STCH_EXPL_Identity2}\\
\left[\beta\right]_{\times}^{2} & =-\beta^{\top}\beta\mathbf{I}_{3}+\beta\beta^{\top},\quad\beta\in\mathbb{R}^{3}\label{eq:SO3STCH_EXPL_Identity3}\\
B\left[\beta\right]_{\times}+ & \left[\beta\right]_{\times}B={\rm Tr}\left\{ B\right\} \left[\beta\right]_{\times}-\left[B\beta\right]_{\times},\nonumber \\
& \quad B=B^{\top}\in\mathbb{R}^{3\times3},\beta\in\mathbb{R}^{3}\label{eq:SO3STCH_EXPL_Identity4}\\
{\rm Tr}\left\{ B\left[\beta\right]_{\times}\right\}  & =0,\quad B=B^{\top}\in\mathbb{R}^{3\times3},\beta\in\mathbb{R}^{3}\label{eq:SO3STCH_EXPL_Identity5}
\end{align}

\section{Problem Formulation \label{sec:SO3STCH_EXP_Problem-Formulation-in}}

\begin{figure}[h]
	\centering{}\includegraphics[scale=0.55]{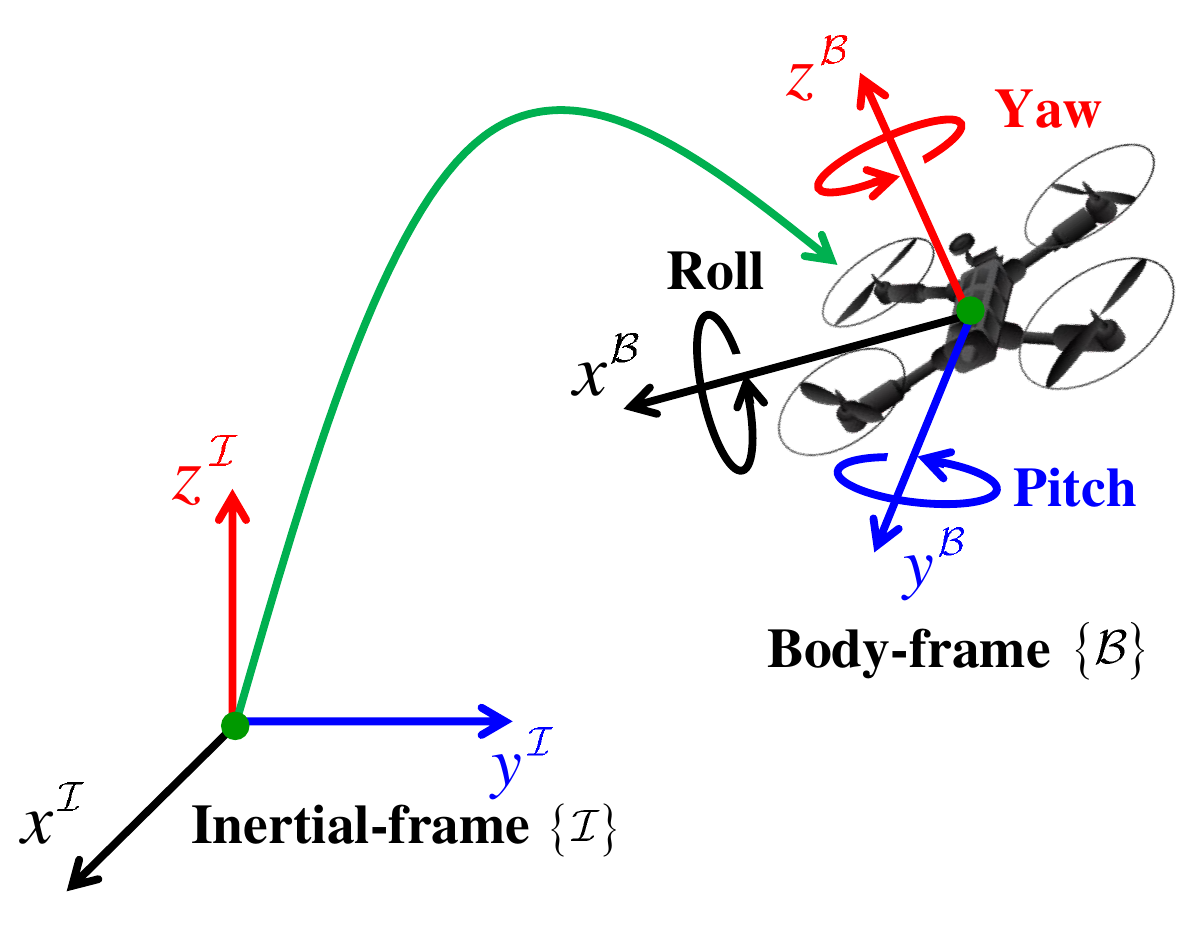}\caption{The relative orientation between body-frame and inertial-frame of
		a rigid-body in 3D space.}
	\label{fig:SO3PPF_1} 
\end{figure}

The attitude can be extracted from $n$-known non-collinear inertial
vectors measured in a coordinate system fixed to the rigid body. Consider
that the superscripts $\mathcal{I}$ and $\mathcal{B}$ refer to the
vectors associated with the inertial-frame and body-frame, respectively.
Let ${\rm v}_{i}^{\mathcal{B}}\in\mathbb{R}^{3}$ be the $i$th measurement
vector in the body-fixed frame for $i=1,2,\ldots,n$. the orientation of the object in the body-frame $\left\{ \mathcal{B}\right\} $
relative to the inertial-frame $\left\{ \mathcal{I}\right\} $ can
be represented by the attitude matrix $R\in\left\{ \mathcal{B}\right\} $
as illustrated in Figure \ref{fig:SO3PPF_1}. 
Let $R\in\mathbb{SO}\left(3\right)$
denote the rotation matrix from body-fixed frame to a given inertial-fixed
frame such that the body-fixed frame vector is defined by
\begin{equation}
{\rm v}_{i}^{\mathcal{B}}=R^{\top}{\rm v}_{i}^{\mathcal{I}}+b_{i}^{\mathcal{B}}+{\rm \omega}_{i}^{\mathcal{B}}\label{eq:SO3STCH_EXPL_Vect_True}
\end{equation}
where ${\rm v}_{i}^{\mathcal{I}}\in\mathbb{R}^{3}$ denotes the inertial-fixed
frame vector while $b_{i}^{\mathcal{B}}$ and ${\rm \omega}_{i}^{\mathcal{B}}$
denote the additive bias and noise components of the associated body-frame
vector, respectively, for all $b_{i}^{\mathcal{B}},{\rm \omega}_{i}^{\mathcal{B}}\in\mathbb{R}^{3}$
and $i=1,2,\ldots,n$. The assumption that $n\ge2$ is necessary for
instantaneous three-dimensional attitude determination. It is common
to employ the normalized values of reference and body-frame vectors
in the process of attitude estimation such as
\begin{equation}
\upsilon_{i}^{\mathcal{I}}=\frac{{\rm v}_{i}^{\mathcal{I}}}{\left\Vert {\rm v}_{i}^{\mathcal{I}}\right\Vert },\hspace{1em}\upsilon_{i}^{\mathcal{B}}=\frac{{\rm v}_{i}^{\mathcal{B}}}{\left\Vert {\rm v}_{i}^{\mathcal{B}}\right\Vert }\label{eq:SO3STCH_EXPL_Vector_norm}
\end{equation}
and the attitude can be defined knowing $\upsilon_{i}^{\mathcal{I}}$
and $\upsilon_{i}^{\mathcal{B}}$. For the sake of simplicity, the
body frame vector (${\rm v}_{i}^{\mathcal{B}}$) is considered to
be noise and bias free in the stability analysis. In the Simulation
Section, on the contrary, noise and bias are present in the measurements.
The true attitude dynamics and the associated Rodriguez vector dynamics
are given in \eqref{eq:SO3STCH_EXPL_R_dynam} and \eqref{eq:SO3STCH_EXPL_Rod_dynam},
respectively, as
\begin{align}
\dot{R} & =R\left[\Omega\right]_{\times}\label{eq:SO3STCH_EXPL_R_dynam}\\
\dot{\rho} & =\frac{1}{2}\left(\mathbf{I}_{3}+\left[\rho\right]_{\times}+\rho\rho^{\top}\right)\Omega\label{eq:SO3STCH_EXPL_Rod_dynam}
\end{align}
where $\Omega\in\mathbb{R}^{3}$ denotes the true value of angular
velocity. Gyroscope or the rate gyros measures the angular velocity
vector in the body-frame relative to the inertial-frame. The measurement
vector of angular velocity is
\begin{equation}
\Omega_{m}=\Omega+b+\omega\label{eq:SO3STCH_EXPL_Angular}
\end{equation}
where $b$ and $\omega$ denote the additive bias and noise components,
respectively, for all $b,\omega\in\mathbb{R}^{3}$. The noise vector
$\omega$ is assumed to be a Gaussian noise vector such that $\mathbb{E}\left[\omega\right]=0$.
The measurement of angular velocity vector is subject to additive
noise and bias, which are characterized by randomness and unknown
behavior, impairing the estimation process of the true attitude dynamics
in \eqref{eq:SO3STCH_EXPL_R_dynam} or \eqref{eq:SO3STCH_EXPL_Rod_dynam}.
As such, \eqref{eq:SO3STCH_EXPL_Angular} is assumed to be excited
by a wide-band of random Gaussian noise process. Derivative of any
Gaussian process yields a Gaussian process allowing the stochastic
attitude dynamics to be written as a function of Brownian motion process
vector \cite{khasminskii1980stochastic,hashim2018SO3Stochastic}
\[
\omega=\mathcal{Q}\frac{d\beta}{dt}
\]
where $\mathcal{Q}\in\mathbb{R}^{3\times3}$ is a non-negative unknown
time-variant diagonal matrix. In addition, each parameter of $\mathcal{Q}$
in the diagonal is bounded with an unknown bound. The properties of
Brownian motion process can be found in \cite{ito1984lectures,deng2001stabilization,hashim2018SO3Stochastic}.
The covariance component associated with noise $\omega$ can be defined
by a diagonal matrix $\mathcal{Q}^{2}=\mathcal{Q}\mathcal{Q}^{\top}$.
Considering the attitude dynamics in \eqref{eq:SO3STCH_EXPL_Rod_dynam}
and replacing $\omega$ by $\mathcal{Q}d\beta/dt$, the stochastic
differential equation in \eqref{eq:SO3STCH_EXPL_Rod_dynam} can be
expressed by
\begin{align}
d\rho= & f\left(\rho,b\right)dt+g\left(\rho\right)\mathcal{Q}d\beta\label{eq:SO3STCH_EXPL_Stoch_Ito}
\end{align}
Similarly, the stochastic dynamics of \eqref{eq:SO3STCH_EXPL_R_dynam}
are
\begin{equation}
dR=R\left[\Omega_{m}-b\right]_{\times}dt-R\left[\mathcal{Q}d\beta\right]_{\times}\label{eq:SO3STCH_EXPL_R_Dyn1}
\end{equation}
where $g\left(\rho\right)=-\frac{1}{2}\left(\mathbf{I}_{3}+\left[\rho\right]_{\times}+\rho\rho^{\top}\right)$
and $f\left(\rho,b\right)=-g\left(\rho\right)\left(\Omega_{m}-b\right)$
with $g:\mathbb{R}^{3}\rightarrow\mathbb{R}^{3\times3}$ and $f:\mathbb{R}^{3}\times\mathbb{R}^{3}\rightarrow\mathbb{R}^{3}$.
$g$ is locally Lipschitz in $\rho$ and $f$ is locally Lipschitz
in $\rho$ and $b$. Accordingly, the dynamic system in \eqref{eq:SO3STCH_EXPL_Stoch_Ito}
has a solution on $t\in\left[t_{0},T\right]\forall t_{0}\leq T<\infty$
in the mean square sense for any $\rho\left(t\right)$ such that $t\neq t_{0}$,
$\rho-\rho_{0}$ is independent of $\left\{ \beta\left(\tau\right),\tau\geq t\right\} ,\forall t\in\left[t_{0},T\right]$
\cite{deng2001stabilization,hashim2018SO3Stochastic}. Aiming to achieve adaptive stabilization
of the unknown time-variant covariance matrix, let us introduce the
following unknown constant
\begin{equation}
\sigma=\left[{\rm max}\left\{ \mathcal{Q}_{1,1}^{2}\right\} ,{\rm max}\left\{ \mathcal{Q}_{2,2}^{2}\right\} ,{\rm max}\left\{ \mathcal{Q}_{3,3}^{2}\right\} \right]^{\top}\label{eq:SO3STCH_EXPL_s_factor}
\end{equation}
where ${\rm max}\left\{ \cdot\right\} $ is the maximum value of the
associated element. From \eqref{eq:SO3STCH_EXPL_Angular}, and \eqref{eq:SO3STCH_EXPL_s_factor},
it can be noticed that $b$ and $\sigma$ are bounded. It is important
to introduce the following Lemma which will be useful in the subsequent
filter derivation. 
\begin{lem}
	\label{Lem:SO3STCH_EXPL_1}Let $R\in\mathbb{SO}\left(3\right)$, $M^{\mathcal{I}}=\left(M^{\mathcal{I}}\right)^{\top}\in\mathbb{R}^{3\times3}$,
	$M^{\mathcal{I}}$ be positive-definite, and ${\rm Tr}\left\{ M^{\mathcal{I}}\right\} =3$.
	Define $\bar{\mathbf{M}}^{\mathcal{I}}={\rm Tr}\left\{ M^{\mathcal{I}}\right\} \mathbf{I}_{3}-M^{\mathcal{I}}$
	and let the minimum singular values of $\bar{\mathbf{M}}^{\mathcal{I}}$
	be $\underline{\lambda}:=\underline{\lambda}\left(\bar{\mathbf{M}}^{\mathcal{I}}\right)$.
	Then, the following holds:
	\begin{align}
	\left\Vert M^{\mathcal{I}}R\right\Vert _{I} & =\frac{1}{2}\frac{\rho^{\top}\bar{\mathbf{M}}^{\mathcal{I}}\rho}{1+\left\Vert \rho\right\Vert ^{2}}\label{eq:SO3STCH_EXPL_lemm1_1}\\
	\Phi\left(M^{\mathcal{I}}R\right) & =\frac{\left(\mathbf{I}_{3}+\left[\rho\right]_{\times}\right)^{\top}\bar{\mathbf{M}}^{\mathcal{I}}}{1+\left\Vert \rho\right\Vert ^{2}}\rho\label{eq:SO3STCH_EXPL_lemm1_2}\\
	\left\Vert M^{\mathcal{I}}R\right\Vert _{I} & \leq\frac{2}{\underline{\lambda}}\frac{||\Phi\left(M^{\mathcal{I}}R\right)||^{2}}{1+{\rm Tr}\left\{ \left(M^{\mathcal{I}}\right)^{-1}M^{\mathcal{I}}R\right\} }\label{eq:SO3STCH_EXPL_lemm1_4}
	\end{align}
	\textbf{Proof. See \nameref{sec:SO3STCH_EXPL_AppendixA}.} 
\end{lem}
\begin{defn}
	\label{def:SO3STCH_EXPL_1}\cite{ji2006adaptive} The Rodriguez vector
	$\rho$ of the stochastic dynamics in \eqref{eq:SO3STCH_EXPL_Stoch_Ito}
	is known to be semi-globally uniformly ultimately bounded (SGUUB)
	if for a compact set $\Lambda\in\mathbb{R}^{3}$ and any $\rho_{0}=\rho\left(t_{0}\right)$,
	there exists a constant $\kappa>0$, and a time constant $T=T\left(\kappa,\rho_{0}\right)$
	such that $\mathbb{E}\left[\left\Vert \rho\right\Vert \right]<\kappa,\forall t>t_{0}+T$. 
\end{defn}
\begin{defn}
	\label{def:SO3STCH_EXPL_2}Consider the stochastic dynamics in \eqref{eq:SO3STCH_EXPL_Stoch_Ito}.
	For a given function $V\left(\rho\right)\in\mathcal{C}^{2}$, the
	differential operator $\mathcal{L}V$ is defined by
	\[
	\mathcal{L}V\left(\rho\right)=V_{\rho}^{\top}f\left(\rho,b\right)+\frac{1}{2}{\rm Tr}\left\{ g\left(\rho\right)\mathcal{Q}^{2}g^{\top}\left(\rho\right)V_{\rho\rho}\right\} 
	\]
	such that $V_{\rho}=\partial V/\partial\rho$, and $V_{\rho\rho}=\partial^{2}V/\partial\rho^{2}$. 
\end{defn}
\begin{lem}
	\label{lem:SO3STCH_EXPL_2} \cite{deng2001stabilization} Let the
	stochastic dynamics in \eqref{eq:SO3STCH_EXPL_Stoch_Ito} be given
	a potential function $V\in\mathcal{C}^{2}$ such that $V:\mathbb{R}^{3}\rightarrow\mathbb{R}^{+}$,
	class $\mathcal{K}_{\infty}$ functions $\varphi_{1}\left(\cdot\right)$
	and $\varphi_{2}\left(\cdot\right)$, constants $c_{1}>0$ and $c_{2}\geq0$,
	and a non-negative function $\boldsymbol{Z}\left(\left\Vert \rho\right\Vert \right)$
	such that
	\begin{equation}
	\varphi_{1}\left(\left\Vert \rho\right\Vert \right)\leq V\left(\rho\right)\leq\varphi_{2}\left(\left\Vert \rho\right\Vert \right)\label{eq:SO3STCH_EXPL_Vfunction_Lyap}
	\end{equation}
	\begin{align}
	\mathcal{L}V\left(\rho\right)= & V_{\rho}^{\top}f\left(\rho,b\right)+\frac{1}{2}{\rm Tr}\left\{ g\left(\rho\right)\mathcal{Q}^{2}g^{\top}\left(\rho\right)V_{\rho\rho}\right\} \nonumber \\
	\leq & -c_{1}\boldsymbol{Z}\left(\left\Vert \rho\right\Vert \right)+c_{2}\label{eq:SO3STCH_EXPL_dVfunction_Lyap}
	\end{align}
	then for $\rho_{0}\in\mathbb{R}^{3}$, there exists almost a unique
	strong solution on $\left[0,\infty\right)$ for the dynamic system
	in \eqref{eq:SO3STCH_EXPL_Stoch_Ito}. The solution $\rho$ is bounded
	in probability such that
	\begin{equation}
	\mathbb{E}\left[V\left(\rho\right)\right]\leq V\left(\rho_{0}\right){\rm exp}\left(-c_{1}t\right)+\frac{c_{2}}{c_{1}}\label{eq:SO3STCH_EXPL_EVfunction_Lyap}
	\end{equation}
	Furthermore, if the inequality in \eqref{eq:SO3STCH_EXPL_EVfunction_Lyap}
	holds, then $\rho$ in \eqref{eq:SO3STCH_EXPL_Stoch_Ito} is SGUUB
	in the mean square. 
\end{lem}
The proof of Lemma \ref{lem:SO3STCH_EXPL_2} can be found in \cite{deng2001stabilization,ji2006adaptive}.
Consider the attitude $R\in\mathbb{SO}\left(3\right)$ and define
the unstable set $\mathcal{U}\subseteq\mathbb{SO}\left(3\right)$
by $\mathcal{U}:=\left\{ \left.R\right|{\rm Tr}\left\{ R\right\} =-1,\boldsymbol{\mathcal{P}}_{a}\left(R\right)=0\right\} $
such that the unstable set $\mathcal{U}$ is forward invariant for
the stochastic dynamics in \eqref{eq:SO3STCH_EXPL_R_dynam} which
implies that $\rho=\infty$. As such, for almost any initial condition
such that $R_{0}\notin\mathcal{U}$ or $\rho_{0}\in\mathbb{R}^{3}$,
one has $-1<{\rm Tr}\left\{ R_{0}\right\} \leq3$ and the trajectory
of $\rho$ converges to the neighborhood of the equilibrium point. 
\begin{lem}
	\label{lem:SO3STCH_EXPL_3} (Young's inequality) Let $y$ and $x$
	be real values such that $y,x\in\mathbb{R}^{n}$. Then, for any $c>0$
	and $d>0$ satisfying $\frac{1}{c}+\frac{1}{d}=1$ with appropriately
	small positive constant $\varepsilon$, the following inequality is
	satisfied
	\begin{align}
	y^{\top}x & \leq\left(1/c\right)\varepsilon^{c}\left\Vert y\right\Vert ^{c}+\left(1/d\right)\varepsilon^{-d}\left\Vert x\right\Vert ^{d}\label{eq:SO3STCH_EXPL_lem_ineq}
	\end{align}
\end{lem}

\section{Nonlinear Stochastic Filter on $\mathbb{SO}\left(3\right)$ \label{sec:SO3STCH_EXP_Stochastic-Complementary-Filters}}

A set of vectorial measurements ${\rm \upsilon}_{i}^{\mathcal{I}}$
and ${\rm \upsilon}_{i}^{\mathcal{B}}$ in \eqref{eq:SO3STCH_EXPL_Vector_norm}
can be employed to reconstruct the uncertain attitude matrix $R_{y}$ such as nonlinear stochastic attitude and pose filters \cite{hashim2018SO3Stochastic,hashim2018SE3Stochastic},
however, obtaining $R_{y}$ might be very computationally complex%
. Therefore, the objective is to propose a nonlinear stochastic attitude
filter which uses a set of vectorial measurements directly without
the need of attitude reconstruction. Consider the error from body-frame
to estimator-frame defined as
\begin{equation}
\tilde{R}=R\hat{R}^{\top}\label{eq:SO3STCH_EXPL_R_error}
\end{equation}
Also, let the error in $b$ and $\sigma$ be given by
\begin{align}
\tilde{b} & =b-\hat{b}\label{eq:SO3STCH_EXPL_b_tilde}\\
\tilde{\sigma} & =\sigma-\hat{\sigma}\label{eq:SO3STCH_EXPL_s_tilde}
\end{align}
Recall $\upsilon_{i}^{\mathcal{I}}$ and $\upsilon_{i}^{\mathcal{B}}$
from \eqref{eq:SO3STCH_EXPL_Vector_norm} for $i=1,\ldots,n$. Define
\begin{align}
M^{\mathcal{I}} & =\left(M^{\mathcal{I}}\right)^{\top}=\sum_{i=1}^{n}s_{i}\upsilon_{i}^{\mathcal{I}}\left(\upsilon_{i}^{\mathcal{I}}\right)^{\top}\nonumber \\
M^{\mathcal{B}} & =\left(M^{\mathcal{B}}\right)^{\top}=\sum_{i=1}^{n}s_{i}\upsilon_{i}^{\mathcal{B}}\left(\upsilon_{i}^{\mathcal{B}}\right)^{\top}=R^{\top}M^{\mathcal{I}}R\label{eq:SO3STCH_EXPL_MB_MI}
\end{align}
with $s_{i}>0$ being the confidence level of the $i$th sensor measurement.
Each of $M^{\mathcal{I}}$ and $M^{\mathcal{B}}$ are symmetric matrices.
Consider $\upsilon_{i}^{\mathcal{I}}$ and $\upsilon_{i}^{\mathcal{B}}$
from \eqref{eq:SO3STCH_EXPL_Vector_norm} for $i=1,\ldots,n$ and
at least two non-collinear vectors available $\left(n\geq2\right)$.
If $n=2$, the third vector is obtained by $\upsilon_{3}^{\mathcal{I}}=\upsilon_{1}^{\mathcal{I}}\times\upsilon_{2}^{\mathcal{I}}$
and $\upsilon_{3}^{\mathcal{B}}=\upsilon_{1}^{\mathcal{B}}\times\upsilon_{2}^{\mathcal{B}}$
which is non-collinear with the other two vectors such that ${\rm rank}\left(M^{\mathcal{I}}\right)={\rm rank}\left(M^{\mathcal{B}}\right)=3$
full rank. Consequently, the three eigenvalues of $M^{\mathcal{I}}$
and $M^{\mathcal{B}}$ are greater than zero. Let $\bar{\mathbf{M}}^{\mathcal{I}}={\rm Tr}\left\{ M^{\mathcal{I}}\right\} \mathbf{I}_{3}-M^{\mathcal{I}}\in\mathbb{R}^{3\times3}$,
provided that ${\rm rank}\left(M^{\mathcal{I}}\right)=3$, the following
statements hold (\cite{bullo2004geometric} page. 553): 
\begin{enumerate}
	\item[i. ] $\bar{\mathbf{M}}^{\mathcal{I}}$ is a symmetric positive-definite
	matrix. 
	\item[ii. ] Define the three eigenvalues of $M^{\mathcal{I}}$ by $\lambda\left(M^{\mathcal{I}}\right)=\left\{ \lambda_{1},\lambda_{2},\lambda_{3}\right\} $,
	then $\lambda\left(\bar{\mathbf{M}}^{\mathcal{I}}\right)=\left\{ \lambda_{3}+\lambda_{2},\lambda_{3}+\lambda_{1},\lambda_{2}+\lambda_{1}\right\} $
	such that the minimum singular value $\underline{\lambda}\left(\bar{\mathbf{M}}^{\mathcal{I}}\right)>0$.
\end{enumerate}
In all the discussion that follows it is assumed that the above statements
hold. Consider $\sum_{i=1}^{n}s_{i}=3$ and define
\begin{equation}
\hat{\upsilon}_{i}^{\mathcal{B}}=\hat{R}^{\top}\upsilon_{i}^{\mathcal{I}}\label{eq:SO3STCH_EXPL_vB_hat}
\end{equation}
From the identity in \eqref{eq:SO3STCH_EXPL_Identity1}, one can find
\begin{align*}
& \sum_{i=1}^{n}\frac{s_{i}}{2}\left[\upsilon_{i}^{\mathcal{B}}\times\hat{\upsilon}_{i}^{\mathcal{B}}\right]_{\times}\\
& =\sum_{i=1}^{n}\frac{s_{i}}{2}\left(\hat{\upsilon}_{i}^{\mathcal{B}}\left(\upsilon_{i}^{\mathcal{B}}\right)^{\top}-\upsilon_{i}^{\mathcal{B}}\left(\hat{\upsilon}_{i}^{\mathcal{B}}\right)^{\top}\right)\\
& =\frac{1}{2}\sum_{i=1}^{n}k_{i}\left(\hat{R}^{\top}\upsilon_{i}^{\mathcal{I}}\left(\upsilon_{i}^{\mathcal{I}}\right)^{\top}R-R^{\top}\upsilon_{i}^{\mathcal{I}}\left(\upsilon_{i}^{\mathcal{I}}\right)^{\top}\hat{R}\right)\\
& =\frac{1}{2}\hat{R}^{\top}\left(M^{\mathcal{I}}\tilde{R}-\tilde{R}^{\top}M^{\mathcal{I}}\right)\hat{R}=\hat{R}^{\top}\boldsymbol{\mathcal{P}}_{a}\left(M^{\mathcal{I}}\tilde{R}\right)\hat{R}
\end{align*}
Hence, the following components can be obtained in terms of vector
measurements which will be used in the proposed filter design
\begin{align}
\Phi\left(M^{\mathcal{I}}\tilde{R}\right)= & \mathbf{vex}\left(\boldsymbol{\mathcal{P}}_{a}\left(M^{\mathcal{I}}\tilde{R}\right)\right)=\hat{R}\sum_{i=1}^{n}\frac{s_{i}}{2}\upsilon_{i}^{\mathcal{B}}\times\hat{\upsilon}_{i}^{\mathcal{B}}\label{eq:SO3STCH_EXPL_VEX_VM}\\
||M^{\mathcal{I}}\tilde{R}||_{I}= & \frac{1}{4}{\rm Tr}\left\{ M^{\mathcal{I}}\left(\mathbf{I}_{3}-\tilde{R}\right)\right\} \nonumber \\
= & \frac{3}{4}-\frac{1}{4}{\rm Tr}\left\{ \hat{R}\sum_{i=1}^{n}\left(s_{i}\hat{\upsilon}_{i}^{\mathcal{B}}\left(\upsilon_{i}^{\mathcal{B}}\right)^{\top}\right)\hat{R}^{\top}\right\} \label{eq:SO3STCH_EXPL_RI_VM}\\
\varUpsilon\left(M^{\mathcal{I}},\tilde{R}\right)= & {\rm Tr}\left\{ \left(\sum_{i=1}^{n}s_{i}\upsilon_{i}^{\mathcal{I}}\left(\upsilon_{i}^{\mathcal{I}}\right)^{\top}\right)^{-1}\right.\nonumber \\
& \hspace{3em}\left.\times\hat{R}\sum_{i=1}^{n}\left(s_{i}\hat{\upsilon}_{i}^{\mathcal{B}}\left(\upsilon_{i}^{\mathcal{B}}\right)^{\top}\right)\hat{R}^{\top}\right\} \label{eq:SO3STCH_EXPL_Trace}
\end{align}
where{\small{}{}{}{}{}$\text{ }\left[\hat{R}\sum_{i=1}^{n}\frac{s_{i}}{2}\upsilon_{i}^{\mathcal{B}}\times\hat{\upsilon}_{i}^{\mathcal{B}}\right]_{\times}=\hat{R}\sum_{i=1}^{n}\frac{s_{i}}{2}\left[\upsilon_{i}^{\mathcal{B}}\times\hat{\upsilon}_{i}^{\mathcal{B}}\right]_{\times}\hat{R}^{\top}$}
as in \eqref{eq:SO3STCH_EXPL_Identity2}. Define $\underline{\lambda}:=\underline{\lambda}\left(\bar{\mathbf{M}}^{\mathcal{I}}\right)$,
$\varUpsilon:=\varUpsilon\left(M^{\mathcal{I}},\tilde{R}\right)$,
and $\Phi:=\Phi\left(M^{\mathcal{I}}\tilde{R}\right)$, and consider
the following nonlinear filter design on $\mathbb{SO}\left(3\right)$
\begin{align}
\dot{\hat{R}}= & \hat{R}\left[\Omega_{m}-\hat{b}\right]_{\times}+\left[W\right]_{\times}\hat{R}\label{eq:SO3STCH_EXPL_Rest_dot_VM}\\
\dot{\hat{b}}= & -\gamma||M^{\mathcal{I}}\tilde{R}||_{I}\hat{R}^{\top}\Phi-\gamma k_{b}\hat{b}\label{eq:SO3STCH_EXPL_b_est_VM}\\
\dot{\hat{\sigma}}= & \frac{\gamma||M^{\mathcal{I}}\tilde{R}||_{I}}{\underline{\text{\ensuremath{\lambda}}}}\frac{{\rm diag}\left(\hat{R}^{\top}\Phi\right)}{1+\varUpsilon}\hat{R}^{\top}\Phi-\gamma k_{\sigma}\hat{\sigma}\label{eq:SO3STCH_EXPL_s_est_VM}\\
W= & \frac{k_{w}}{\varepsilon\underline{\lambda}}\frac{\left(1+\varUpsilon\right)^{2}\underline{\lambda}^{2}+1}{1+\varUpsilon}\Phi+\frac{1}{\underline{\text{\ensuremath{\lambda}}}}\frac{\hat{R}{\rm diag}\left(\hat{R}^{\top}\Phi\right)}{1+\varUpsilon}\hat{\sigma}\label{eq:SO3STCH_EXPL_W_VM}
\end{align}
where $\Phi\left(M^{\mathcal{I}}\tilde{R}\right)$, $\left\Vert M^{\mathcal{I}}\tilde{R}\right\Vert _{I}$,
and $\varUpsilon\left(M^{\mathcal{I}},\tilde{R}\right)$ are defined
in \eqref{eq:SO3STCH_EXPL_VEX_VM}, \eqref{eq:SO3STCH_EXPL_RI_VM},
and \eqref{eq:SO3STCH_EXPL_Trace} in terms of vectorial measurements,
respectively, ${\rm diag}\left(\cdot\right)$ is a diagonal of the
associated component, $k_{w}$, $k_{b},$ $k_{\sigma}$, and $\gamma$
are positive constants, and $\hat{b}$ and $\hat{\sigma}$ are the
estimate of $b$ and $\sigma$, respectively. 
\begin{thm}
	\textbf{\label{thm:SO3STCH_EXPL_2}} Consider the observer in \eqref{eq:SO3STCH_EXPL_Rest_dot_VM},
	\eqref{eq:SO3STCH_EXPL_b_est_VM}, \eqref{eq:SO3STCH_EXPL_s_est_VM}
	and \eqref{eq:SO3STCH_EXPL_W_VM} coupled with angular velocity measurements
	in \eqref{eq:SO3STCH_EXPL_Angular} and the normalized vectors in
	\eqref{eq:SO3STCH_EXPL_Vector_norm}. Assume that two or more body-frame
	non-collinear vectors are available for measurements such that $M^{\mathcal{I}}$
	in \eqref{eq:SO3STCH_EXPL_MB_MI} is nonsingular. Then, for angular
	velocity measurements contaminated with noise and $\tilde{\rho}\in\mathbb{R}^{3}$,
	$\tilde{\rho}$, $\tilde{b}$ and $\tilde{\sigma}$ are regulated
	to an arbitrarily small neighborhood of the origin in probability;
	and $\left[\tilde{\rho}^{\top},\tilde{b}^{\top},\tilde{\sigma}^{\top}\right]^{\top}$
	is SGUUB in mean square. 
\end{thm}
\textbf{Proof: }Let the error in attitude be $\tilde{R}=R\hat{R}^{\top}$
as given in \eqref{eq:SO3STCH_EXPL_R_error} and consider \eqref{eq:SO3STCH_EXPL_b_tilde}
and \eqref{eq:SO3STCH_EXPL_s_tilde}. In view of \eqref{eq:SO3STCH_EXPL_R_dynam}
and \eqref{eq:SO3STCH_EXPL_Rest_dot_VM}, the derivative of attitude
error in incremental form becomes
\begin{align}
d\tilde{R}= & -R\left[\Omega_{m}-\hat{b}\right]_{\times}\hat{R}^{\top}dt-R\hat{R}^{\top}\left[W\right]_{\times}dt\nonumber \\
& +R\left[\Omega_{m}-b-\mathcal{Q}\frac{d\beta}{dt}\right]_{\times}\hat{R}^{\top}dt\nonumber \\
= & -R\left[\tilde{\sigma}\right]_{\times}\hat{R}^{\top}dt-R\hat{R}^{\top}\left[W\right]_{\times}dt-R\left[\mathcal{Q}d\beta\right]_{\times}\hat{R}^{\top}\nonumber \\
= & -\tilde{R}\left[\hat{R}\tilde{b}+W\right]_{\times}dt-\tilde{R}\left[\hat{R}\mathcal{Q}d\beta\right]_{\times}\label{eq:SO3STCH_EXPL_Rtilde_dot}
\end{align}
where $\left[\hat{R}\tilde{\sigma}\right]_{\times}=\hat{R}\left[\tilde{\sigma}\right]_{\times}\hat{R}^{\top}$
as in \eqref{eq:SO3STCH_EXPL_Identity2}. Recalling the extraction
of Rodriguez vector dynamics from \eqref{eq:SO3STCH_EXPL_R_Dyn1}
to \eqref{eq:SO3STCH_EXPL_Stoch_Ito}, the Rodriguez error vector
dynamic in \eqref{eq:SO3STCH_EXPL_Rtilde_dot} can be expressed as
\begin{align}
d\tilde{\rho}= & \tilde{f}dt+\tilde{g}\hat{R}\mathcal{Q}d\beta\label{eq:SO3STCH_EXPL_Rod_tilde_dot_Ito}
\end{align}
where $\tilde{\rho}$ is a Rodriguez error vector associated with
$\tilde{R}$, $\tilde{g}=-\frac{1}{2}\left(\mathbf{I}_{3}+\left[\tilde{\rho}\right]_{\times}+\tilde{\rho}\tilde{\rho}^{\top}\right)$,
and $\tilde{f}=\tilde{g}\left(\hat{R}\tilde{b}+W\right)$. 
\begin{rem}
	\label{rem:The-traditional-potential} From literature, one of the
	traditional potential functions of the adaptive filter is similar
	to \cite{mahony2008nonlinear,crassidis2007survey,zlotnik2017nonlinear}
	\begin{equation}
	V\left(\tilde{R},\tilde{b}\right)=\frac{1}{4}{\rm Tr}\left\{ M^{\mathcal{I}}\left(\mathbf{I}_{3}-\tilde{R}\right)\right\} +\frac{1}{2\gamma}\tilde{b}^{\top}\tilde{b}\label{eq:SO3STCH_EXPL_LyapV-R}
	\end{equation}
	Given \eqref{eq:SO3STCH_EXPL_lemm1_1}, the expression in \eqref{eq:SO3STCH_EXPL_LyapV-R}
	is equivalent to \eqref{eq:SO3STCH_EXPL_LyapV-Trad} in Rodriquez
	vector form
	\begin{align}
	V\left(\tilde{\rho},\tilde{b}\right)= & \frac{1}{2}\frac{\tilde{\rho}^{\top}\bar{\mathbf{M}}^{\mathcal{I}}\tilde{\rho}}{1+\left\Vert \tilde{\rho}\right\Vert ^{2}}+\frac{1}{2\gamma}\tilde{b}^{\top}\tilde{b}\label{eq:SO3STCH_EXPL_LyapV-Trad}
	\end{align}
	The weakness of the potential function in \eqref{eq:SO3STCH_EXPL_LyapV-Trad}
	is that the trace component of the operator $\mathcal{L}V$ in Definition
	\ref{def:SO3STCH_EXPL_2} for the dynamic system in \eqref{eq:SO3STCH_EXPL_Stoch_Ito}
	at $\tilde{\rho}=0$ can be evaluated by
	\[
	\left.\frac{1}{2}{\rm Tr}\left\{ \hat{R}^{\top}\tilde{g}^{\top}V_{\tilde{\rho}\tilde{\rho}}\tilde{g}\hat{R}\mathcal{Q}^{2}\right\} \right|_{\tilde{\rho}=0}=\frac{1}{8}{\rm Tr}\left\{ \hat{R}^{\top}\bar{\mathbf{M}}^{\mathcal{I}}\hat{R}\mathcal{Q}^{2}\right\} 
	\]
	such that the significant impact of $\mathcal{Q}^{2}$ cannot be lessened. 
\end{rem}
Therefore, consider the following smooth attitude potential function
\begin{align}
V\left(\tilde{\rho},\tilde{b},\tilde{\sigma}\right)= & \frac{1}{4}\left(\frac{\tilde{\rho}^{\top}\bar{\mathbf{M}}^{\mathcal{I}}\tilde{\rho}}{1+\left\Vert \tilde{\rho}\right\Vert ^{2}}\right)^{2}+\frac{1}{2\gamma}\tilde{b}^{\top}\tilde{b}+\frac{1}{2\gamma}\tilde{\sigma}^{\top}\tilde{\sigma}\label{eq:SO3STCH_EXPL_LyapV}
\end{align}
For $V:=V\left(\tilde{\rho},\tilde{b},\tilde{\sigma}\right)$, the
differential operator $\mathcal{L}V$ in Definition \ref{def:SO3STCH_EXPL_2}
can be written a{\small{}s
	\begin{equation}
	\mathcal{L}V=V_{\tilde{\rho}}^{\top}\tilde{f}+\frac{1}{2}{\rm Tr}\left\{ \hat{R}^{\top}\tilde{g}^{\top}V_{\tilde{\rho}\tilde{\rho}}\tilde{g}\hat{R}\mathcal{Q}^{2}\right\} -\frac{1}{\gamma}\tilde{b}^{\top}\dot{\hat{b}}-\frac{1}{\gamma}\tilde{\sigma}^{\top}\dot{\hat{\sigma}}\label{eq:SO3STCH_EXPL_LyapLV}
	\end{equation}
}Hence, the first and second partial derivatives of \eqref{eq:SO3STCH_EXPL_LyapV}
can be defined respectively, as follows
\begin{align}
V_{\mathcal{\tilde{\rho}}}= & \frac{\tilde{\rho}^{\top}\bar{\mathbf{M}}^{\mathcal{I}}\tilde{\rho}}{\left(1+\left\Vert \tilde{\rho}\right\Vert ^{2}\right)^{3}}\left(\left(1+\left\Vert \tilde{\rho}\right\Vert ^{2}\right)\mathbf{I}_{3}-\tilde{\rho}\tilde{\rho}^{\top}\right)\bar{\mathbf{M}}^{\mathcal{I}}\tilde{\rho}\label{eq:SO3STCH_EXPL_LyapVv}\\
V_{\tilde{\rho}\tilde{\rho}}= & \frac{\tilde{\rho}^{\top}\bar{\mathbf{M}}^{\mathcal{I}}\tilde{\rho}}{\left(1+\left\Vert \tilde{\rho}\right\Vert ^{2}\right)^{2}}\bar{\mathbf{M}}^{\mathcal{I}}+2\frac{\bar{\mathbf{M}}^{\mathcal{I}}\tilde{\rho}\tilde{\rho}^{\top}\bar{\mathbf{M}}^{\mathcal{I}}}{\left(1+\left\Vert \tilde{\rho}\right\Vert ^{2}\right)^{2}}\nonumber \\
& -4\frac{\tilde{\rho}^{\top}\bar{\mathbf{M}}^{\mathcal{I}}\tilde{\rho}}{\left(1+\left\Vert \tilde{\rho}\right\Vert ^{2}\right)^{3}}\left(\bar{\mathbf{M}}^{\mathcal{I}}\tilde{\rho}\tilde{\rho}^{\top}+\tilde{\rho}\tilde{\rho}^{\top}\bar{\mathbf{M}}^{\mathcal{I}}\right)\nonumber \\
& +\frac{\left(\tilde{\rho}^{\top}\bar{\mathbf{M}}^{\mathcal{I}}\tilde{\rho}\right)^{2}}{\left(1+\left\Vert \tilde{\rho}\right\Vert ^{2}\right)^{4}}\left(6\tilde{\rho}\tilde{\rho}^{\top}-\left(1+\left\Vert \tilde{\rho}\right\Vert ^{2}\right)\mathbf{I}_{3}\right)\label{eq:SO3STCH_EXPL_LyapVvv}
\end{align}
from \eqref{eq:SO3STCH_EXPL_Rod_tilde_dot_Ito} and \eqref{eq:SO3STCH_EXPL_LyapVv},
the first part of \eqref{eq:SO3STCH_EXPL_LyapLV} can be defined by
{\small{}{}{}{}{}
	\begin{align}
	V_{\mathcal{\tilde{\rho}}}^{\top}\tilde{f} & =-\frac{1}{2}\frac{\tilde{\rho}^{\top}\bar{\mathbf{M}}^{\mathcal{I}}\tilde{\rho}\tilde{\rho}^{\top}\bar{\mathbf{M}}^{\mathcal{I}}}{\left(1+\left\Vert \tilde{\rho}\right\Vert ^{2}\right)^{2}}\left(\mathbf{I}_{3}+\left[\tilde{\rho}\right]_{\times}\right)\left(\hat{R}\tilde{b}+W\right)dt\nonumber \\
	& =-||M^{\mathcal{I}}\tilde{R}||_{I}\Phi\left(M^{\mathcal{I}}\tilde{R}\right)^{\top}\left(\hat{R}\tilde{b}+W\right)dt\label{eq:SO3STCH_EXPL_LyapLV_Part1}
	\end{align}
}where $||M^{\mathcal{I}}\tilde{R}||_{I}$ and $\Phi\left(M^{\mathcal{I}}\tilde{R}\right)$
are defined in \eqref{eq:SO3STCH_EXPL_lemm1_1} and \eqref{eq:SO3STCH_EXPL_lemm1_2},
respectively. From \eqref{eq:SO3STCH_EXPL_Rod_tilde_dot_Ito} and
\eqref{eq:SO3STCH_EXPL_LyapVvv}, the second part of \eqref{eq:SO3STCH_EXPL_LyapLV}
can be obtained by{\small{}
	\begin{align}
	& \frac{1}{2}{\rm Tr}\left\{ \hat{R}^{\top}\tilde{g}^{\top}V_{\tilde{\rho}\tilde{\rho}}\tilde{g}\hat{R}\mathcal{Q}^{2}\right\} =-\frac{1}{4}{\rm Tr}\left\{ \frac{1}{4}\left(\frac{\tilde{\rho}^{\top}\bar{\mathbf{M}}^{\mathcal{I}}\tilde{\rho}}{1+\left\Vert \tilde{\rho}\right\Vert ^{2}}\right)^{2}\hat{R}\mathcal{Q}^{2}\hat{R}^{\top}\right\} \nonumber \\
	& +\frac{1}{8}\frac{\tilde{\rho}^{\top}\bar{\mathbf{M}}^{\mathcal{I}}\tilde{\rho}}{\left(1+\left\Vert \tilde{\rho}\right\Vert ^{2}\right)^{2}}{\rm Tr}\left\{ \left(\mathbf{I}_{3}+\left[\tilde{\rho}\right]_{\times}\right)^{\top}\bar{\mathbf{M}}^{\mathcal{I}}\left(\mathbf{I}_{3}+\left[\tilde{\rho}\right]_{\times}\right)\hat{R}\mathcal{Q}^{2}\hat{R}^{\top}\right.\nonumber \\
	& \hspace{1em}\left.-\left(\tilde{\rho}\tilde{\rho}^{\top}\bar{\mathbf{M}}^{\mathcal{I}}\left(\mathbf{I}_{3}+\left[\tilde{\rho}\right]_{\times}\right)+\left(\mathbf{I}_{3}+\left[\tilde{\rho}\right]_{\times}\right)^{\top}\bar{\mathbf{M}}^{\mathcal{I}}\tilde{\rho}\tilde{\rho}^{\top}\right)\hat{R}\mathcal{Q}^{2}\hat{R}^{\top}\right\} \nonumber \\
	& +{\rm Tr}\left\{ \frac{\left(\mathbf{I}_{3}+\left[\tilde{\rho}\right]_{\times}\right)^{\top}\bar{\mathbf{M}}^{\mathcal{I}}\tilde{\rho}\tilde{\rho}^{\top}\bar{\mathbf{M}}^{\mathcal{I}}\left(\mathbf{I}_{3}+\left[\tilde{\rho}\right]_{\times}\right)}{4\left(1+\left\Vert \tilde{\rho}\right\Vert ^{2}\right)^{2}}\hat{R}\mathcal{Q}^{2}\hat{R}^{\top}\right\} \label{eq:SO3STCH_EXPL_LyapLV_Part2_1}
	\end{align}
}from \eqref{eq:SO3STCH_EXPL_lemm1_1} and \eqref{eq:SO3STCH_EXPL_lemm1_2},
one has{\small{}
	\begin{align}
	& \frac{1}{2}{\rm Tr}\left\{ \hat{R}^{\top}\tilde{g}^{\top}V_{\tilde{\rho}\tilde{\rho}}\tilde{g}\hat{R}\mathcal{Q}^{2}\right\} =\nonumber \\
	& -\frac{1}{4}{\rm Tr}\left\{ ||M^{\mathcal{I}}\tilde{R}||_{I}\left(\frac{\left(\mathbf{I}_{3}+\left[\tilde{\rho}\right]_{\times}\right)^{\top}M^{\mathcal{I}}\left(\mathbf{I}_{3}+\left[\tilde{\rho}\right]_{\times}\right)}{1+\left\Vert \tilde{\rho}\right\Vert ^{2}}\right.\right.\nonumber \\
	& \hspace{8em}\left.\left.+||M^{\mathcal{I}}\tilde{R}||_{I}\mathbf{I}_{3}\right)\hat{R}\mathcal{Q}^{2}\hat{R}^{\top}\right\} \nonumber \\
	& +\frac{1}{4}{\rm Tr}\left\{ \left(\Phi\Phi^{\top}+||M^{\mathcal{I}}\tilde{R}||_{I}\left(3\mathbf{I}_{3}-2\Phi\tilde{\rho}^{\top}\right)\right)\hat{R}\mathcal{Q}^{2}\hat{R}^{\top}\right\} \label{eq:SO3STCH_EXPL_LyapLV_Part2_2}
	\end{align}
}where the first part of \eqref{eq:SO3STCH_EXPL_LyapLV_Part2_2} is
negative for all $\tilde{\rho}\neq0$ and $\mathcal{Q}^{2}\neq0$,
and $\Phi:=\Phi\left(M^{\mathcal{I}}\tilde{R}\right)$. From \eqref{eq:SO3STCH_EXPL_TR2}
and \eqref{eq:SO3STCH_EXPL_VEX_Pa}, one can easily find that for
$\varUpsilon:=\varUpsilon\left(M^{\mathcal{I}},\tilde{R}\right)$
\begin{equation}
1+||\tilde{\rho}||^{2}=\frac{1}{1-||\tilde{R}||_{I}}=\frac{4}{1+\varUpsilon}\label{eq:SO3STCH_EXPL_LyapLV_Part2_3}
\end{equation}
Accordingly, from \eqref{eq:SO3STCH_EXPL_VEX_Pa}, $\Phi\left(\tilde{R}\right)=2\tilde{\rho}/\left(1+||\tilde{\rho}||^{2}\right)$,
and from \eqref{eq:SO3STCH_EXPL_lemm1_2}, $\Phi\left(M^{\mathcal{I}}\tilde{R}\right)=\left(\mathbf{I}_{3}+\left[\tilde{\rho}\right]_{\times}\right)^{\top}\bar{\mathbf{M}}^{\mathcal{I}}\tilde{\rho}/\left(1+||\tilde{\rho}||^{2}\right)$.
In addition to the result in \eqref{eq:SO3STCH_EXPL_LyapLV_Part2_3},
one has
\begin{equation}
\underline{\text{\ensuremath{\lambda}}}\Phi^{\top}\hat{R}\mathcal{Q}^{2}\hat{R}^{\top}\tilde{\rho}\leq2\frac{\Phi^{\top}\hat{R}{\rm diag}\left(\hat{R}^{\top}\Phi\right)}{1+\varUpsilon}\sigma\label{eq:SO3STCH_EXPL_LyapLV_Part2_4}
\end{equation}
Define $q=\left[\mathcal{Q}_{1,1},\mathcal{Q}_{2,2},\mathcal{Q}_{3,3}\right]^{\top}$,
as ${\rm Tr}\left\{ \hat{R}\mathcal{Q}^{2}\hat{R}^{\top}\right\} ={\rm Tr}\left\{ \mathcal{Q}^{2}\right\} $,
thereby, the following inequality holds
\begin{equation}
{\rm Tr}\left\{ \Phi\Phi^{\top}\hat{R}\mathcal{Q}^{2}\hat{R}^{\top}\right\} \leq\left\Vert q\right\Vert ^{2}\left\Vert \Phi\right\Vert ^{2}\label{eq:SO3STCH_EXPL_LyapLV_Part2_5}
\end{equation}
Let us combine the results in \eqref{eq:SO3STCH_EXPL_LyapLV_Part2_4}
and \eqref{eq:SO3STCH_EXPL_LyapLV_Part2_5} with \eqref{eq:SO3STCH_EXPL_LyapLV_Part2_2}.
Next, we express the differential operator in \eqref{eq:SO3STCH_EXPL_LyapLV}
in its complete for{\small{}m
	\begin{align}
	\mathcal{L}V\leq & -||M^{\mathcal{I}}\tilde{R}||_{I}\Phi^{\top}\left(\hat{R}\tilde{b}+W\right)\nonumber \\
	& -\frac{1}{4}{\rm Tr}\left\{ ||M^{\mathcal{I}}\tilde{R}||_{I}\left(\frac{\left(\mathbf{I}_{3}+\left[\tilde{\rho}\right]_{\times}\right)^{\top}M^{\mathcal{I}}\left(\mathbf{I}_{3}+\left[\tilde{\rho}\right]_{\times}\right)}{1+\left\Vert \tilde{\rho}\right\Vert ^{2}}\right.\right.\nonumber \\
	& \hspace{7em}\left.\left.+||M^{\mathcal{I}}\tilde{R}||_{I}\mathbf{I}_{3}\right)\hat{R}\mathcal{Q}^{2}\hat{R}^{\top}\right\} \nonumber \\
	& +\frac{1}{4}{\rm Tr}\left\{ \left(\left\Vert \Phi\right\Vert ^{2}+3||M^{\mathcal{I}}\tilde{R}||_{I}\right)\left\Vert q\right\Vert ^{2}\right\} \nonumber \\
	& +\frac{||M^{\mathcal{I}}\tilde{R}||_{I}}{\underline{\text{\ensuremath{\lambda}}}}\Phi^{\top}\hat{R}\frac{{\rm diag}\left(\hat{R}^{\top}\Phi\right)}{1+\varUpsilon}\sigma-\frac{1}{\gamma}\tilde{b}^{\top}\dot{\hat{b}}-\frac{1}{\gamma}\tilde{\sigma}^{\top}\dot{\hat{\sigma}}\label{eq:SO3STCH_EXPL_LyapLV_Part2_6}
	\end{align}
}Considering \eqref{eq:SO3STCH_EXPL_lem_ineq} in Lemma \ref{lem:SO3STCH_EXPL_3},
one obtains
\begin{align}
\left\Vert q\right\Vert ^{2}\left\Vert \Phi\right\Vert ^{2} & \leq\frac{\varepsilon}{2}\left\Vert q\right\Vert ^{4}+\frac{1}{2\varepsilon}\left\Vert \Phi\right\Vert ^{4}\nonumber \\
\left\Vert q\right\Vert ^{2}||M^{\mathcal{I}}\tilde{R}||_{I} & \leq\frac{\varepsilon}{2}\left\Vert q\right\Vert ^{4}+\frac{1}{2\varepsilon}||M^{\mathcal{I}}\tilde{R}||_{I}^{2}\label{eq:SO3STCH_EXPL_LyapLV_Part2_7}
\end{align}
since the second term in \eqref{eq:SO3STCH_EXPL_LyapLV_Part2_6} is
negative semi-definite, we combine \eqref{eq:SO3STCH_EXPL_LyapLV_Part2_7}
with \eqref{eq:SO3STCH_EXPL_LyapLV_Part2_6}. Disregarding the second
term in \eqref{eq:SO3STCH_EXPL_LyapLV_Part2_6} and consider the inequality
in \eqref{eq:SO3STCH_EXPL_lemm1_4} such {\small{}that
	\begin{align}
	& \mathcal{L}V\leq-||M^{\mathcal{I}}\tilde{R}||_{I}\Phi^{\top}\left(\hat{R}\tilde{b}+W\right)\nonumber \\
	& \:+\frac{1}{\underline{\text{\ensuremath{\lambda}}}}||M^{\mathcal{I}}\tilde{R}||_{I}\Phi^{\top}\left(\frac{1}{2\varepsilon}\frac{\left(1+\varUpsilon\right)^{2}\underline{\lambda}^{2}+1}{1+\varUpsilon}\Phi+\frac{\hat{R}{\rm diag}\left(\hat{R}^{\top}\Phi\right)}{1+\varUpsilon}\sigma\right)\nonumber \\
	& \:-\frac{1}{\gamma}\tilde{b}^{\top}\dot{\hat{b}}-\frac{1}{\gamma}\tilde{\sigma}^{\top}\dot{\hat{\sigma}}+\frac{\varepsilon}{2}\bar{\sigma}^{2}\label{eq:SO3STCH_EXPL_LyapLV_Part2_8}
	\end{align}
}where $\bar{\sigma}=\sum_{i=1}^{3}\sigma_{i}\geq\left\Vert q\right\Vert ^{2}$.
With direct substitution of $\dot{\hat{b}}$, $\dot{\hat{\sigma}}$
and $W$ from \eqref{eq:SO3STCH_EXPL_b_est_VM}, \eqref{eq:SO3STCH_EXPL_s_est_VM},
and \eqref{eq:SO3STCH_EXPL_W_VM}, respectively, one finds
\begin{align}
\mathcal{L}V\leq & -\frac{2k_{w}-1}{2\varepsilon}\left(\underline{\lambda}^{2}\left(1+\varUpsilon\right)^{2}+1\right)||M^{\mathcal{I}}\tilde{R}||_{I}^{2}-k_{b}||\tilde{b}||^{2}\nonumber \\
& -k_{\sigma}\left\Vert \tilde{\sigma}\right\Vert ^{2}+k_{b}\tilde{b}^{\top}b+k_{\sigma}\tilde{\sigma}^{\top}\sigma+\frac{\varepsilon}{2}\bar{\sigma}^{2}\label{eq:SO3STCH_EXPL_LyapLV_Part2_9}
\end{align}
According to \eqref{eq:SO3STCH_EXPL_lem_ineq} in Lemma \ref{lem:SO3STCH_EXPL_3},
one has
\begin{align*}
\tilde{b}^{\top}b & \leq\frac{1}{2}||\tilde{b}||^{2}+\frac{1}{2}\left\Vert b\right\Vert ^{2}\\
\tilde{\sigma}^{\top}\sigma & \leq\frac{1}{2}\left\Vert \tilde{\sigma}\right\Vert ^{2}+\frac{1}{2}\left\Vert \sigma\right\Vert ^{2}
\end{align*}
Thus, the differential operator in \eqref{eq:SO3STCH_EXPL_LyapLV_Part2_9}
becomes
\begin{align}
\mathcal{L}V\leq & -\frac{2k_{w}-1}{2\varepsilon}||M^{\mathcal{I}}\tilde{R}||_{I}^{2}-\frac{k_{b}}{2}||\tilde{b}||^{2}-\frac{k_{\sigma}}{2}\left\Vert \tilde{\sigma}\right\Vert ^{2}\nonumber \\
& +\frac{1}{2}\left(k_{\sigma}+\varepsilon\right)\bar{\sigma}^{2}+\frac{1}{2}k_{b}\left\Vert b\right\Vert ^{2}\label{eq:SO3STCH_EXPL_LyapLV_Part2_10}
\end{align}
Define%
\begin{align*}
c_{2} & =\frac{1}{2}\left(k_{\sigma}+\varepsilon\right)\bar{\sigma}^{2}+\frac{1}{2}k_{b}\left\Vert b\right\Vert ^{2}\in\mathbb{R}\\
\tilde{X} & =\left[\frac{1}{2}\frac{\tilde{\rho}^{\top}\bar{\mathbf{M}}^{\mathcal{I}}\tilde{\rho}}{1+\left\Vert \tilde{\rho}\right\Vert ^{2}},\frac{1}{\sqrt{2\gamma}}\tilde{b}^{\top},\frac{1}{\sqrt{2\gamma}}\tilde{\sigma}^{\top}\right]^{\top}\in\mathbb{R}^{7}\\
\mathcal{H} & =\left[\begin{array}{ccc}
\frac{2k_{w}-1}{2\varepsilon} & \underline{\mathbf{0}}_{3}^{\top} & \underline{\mathbf{0}}_{3}^{\top}\\
\underline{\mathbf{0}}_{3} & \gamma k_{b}\mathbf{I}_{3} & \mathbf{0}_{3\times3}\\
\underline{\mathbf{0}}_{3} & \mathbf{0}_{3\times3} & \gamma k_{\sigma}\mathbf{I}_{3}
\end{array}\right]\in\mathbb{R}^{7\times7}
\end{align*}
as such, the differential operator in \eqref{eq:SO3STCH_EXPL_LyapLV_Part2_10}
becomes
\begin{align}
\mathcal{L}V & \leq-\tilde{X}^{\top}\mathcal{H}\tilde{X}+c_{2}\leq-\underline{\lambda}\left(\mathcal{H}\right)V+c_{2}\label{eq:SO3STCH_EXPL_LyapLV_Final}
\end{align}
where $\underline{\lambda}\left(\cdot\right)$ is the minimum singular
value of a matrix. Hence, from \eqref{eq:SO3STCH_EXPL_LyapLV_Final},
one has
\begin{equation}
d\left(\mathbb{E}\left[V\right]\right)/dt=\mathbb{E}\left[\mathcal{L}V\right]\leq-\underline{\lambda}\left(\mathcal{H}\right)V+c_{2}\label{eq:SO3STCH_EXPL_dV_Exp_ito}
\end{equation}
According to Lemma \eqref{lem:SO3STCH_EXPL_2}, the inequality in
\eqref{eq:SO3STCH_EXPL_dV_Exp_ito} means{\small{}
	\begin{equation}
	0\leq\mathbb{E}\left[V\left(t\right)\right]\leq V\left(0\right){\rm exp}\left(-\underline{\lambda}\left(\mathcal{H}\right)t\right)+\frac{c_{2}}{\underline{\lambda}\left(\mathcal{H}\right)},\,\forall t\geq0\label{eq:SO3STCH_EXPL_V_Exp_ito}
	\end{equation}
}The inequality in \eqref{eq:SO3STCH_EXPL_V_Exp_ito} means that $\mathbb{E}\left[V\left(t\right)\right]$
is ultimately bounded by $c_{2}/\underline{\lambda}\left(\mathcal{H}\right)$%
. Let $\tilde{Y}=[\tilde{\rho}^{\top},\tilde{b}^{\top},\tilde{\sigma}^{\top}]^{\top}$,
hence, $\tilde{Y}$ is SGUUB in the mean square.%
{} For $\tilde{Y}_{0}\in\mathbb{R}^{9}$, the trajectory of $\tilde{Y}$
steers to the neighborhood of the origin and $c_{2}/\underline{\lambda}\left(\mathcal{H}\right)$
being the ultimate upper bound of the neighborhood.

\section{Simulation \label{sec:SO3STCH_EXP_Simulation}}

Let $R$ be expressed by the dynamics in \eqref{eq:SO3STCH_EXPL_R_dynam}
with {\small{}$\Omega=\left[{\rm sin}\left(0.7t\right),0.7{\rm sin}\left(0.5t+\pi\right),0.5{\rm sin}\left(0.3t+\frac{\pi}{3}\right)\right]^{\top}{\rm rad/sec}$}
and initial attitude $R\left(0\right)=\mathbf{I}_{3}$. The true angular velocity is considered to be corrupted by a wide-band
of random noise process $\omega$ with standard deviation (STD) being
$0.2\left({\rm rad/sec}\right)$ and zero mean, and bias $b=0.2\left[1,-1,1\right]^{\top}$.
Consider two non-collinear inertial-frame measurements being given
by ${\rm v}_{1}^{\mathcal{I}}=\frac{1}{\sqrt{3}}\left[1,-1,1\right]^{\top}$
and ${\rm v}_{2}^{\mathcal{I}}=\left[0,0,1\right]^{\top}$ and their
body-frame measurements being given by ${\rm v}_{i}^{\mathcal{B}}=R^{\top}{\rm v}_{i}^{\mathcal{I}}+{\rm b}_{i}^{\mathcal{B}}+\omega_{i}^{\mathcal{B}}$
where $\omega_{1}^{\mathcal{B}}$ and $\omega_{2}^{\mathcal{B}}$
are Gaussian noise process vectors with ${\rm STD}=0.2$ and zero
mean and the associated bias components ${\rm b}_{1}^{\mathcal{B}}=0.1\left[-1,1,0.5\right]^{\top}$
and ${\rm b}_{2}^{\mathcal{B}}=0.1\left[0,0,1\right]^{\top}$. The
third vector is obtained by the cross product.

$\hat{R}\left(0\right)$ is given by angle-axis parameterization in
\eqref{eq:SO3STCH_EXPL_att_ang} as $\hat{R}\left(0\right)=\mathcal{R}_{\alpha}\left(\alpha,u/\left\Vert u\right\Vert \right)$
with $\alpha=179\left({\rm deg}\right)$ and $u$= $\left[1,5,3\right]^{\top}$
such that $\tilde{R}$ approaches the unstable equilibria $||\tilde{R}||_{I}\approx0.9999$
\[
R\left(0\right)=\mathbf{I}_{3},\hspace{1em}\hat{R}\left(0\right)=\left[\begin{array}{ccc}
-0.9427 & 0.2768 & 0.1862\\
0.2945 & 0.4286 & 0.8541\\
0.1567 & 0.8600 & -0.4856
\end{array}\right]
\]
Initial estimates are selected as $\hat{b}\left(0\right)=\underline{\mathbf{0}}_{3}$,
$\hat{\sigma}\left(0\right)=\underline{\mathbf{0}}_{3}$, and design
parameters are as follows: $\gamma=1$, $k_{b}=0.5$, $k_{\sigma}=0.5$,
$k_{w}=5$ , and $\varepsilon=0.5$.

\begin{figure}[h!]
	\centering{}\includegraphics[scale=0.26]{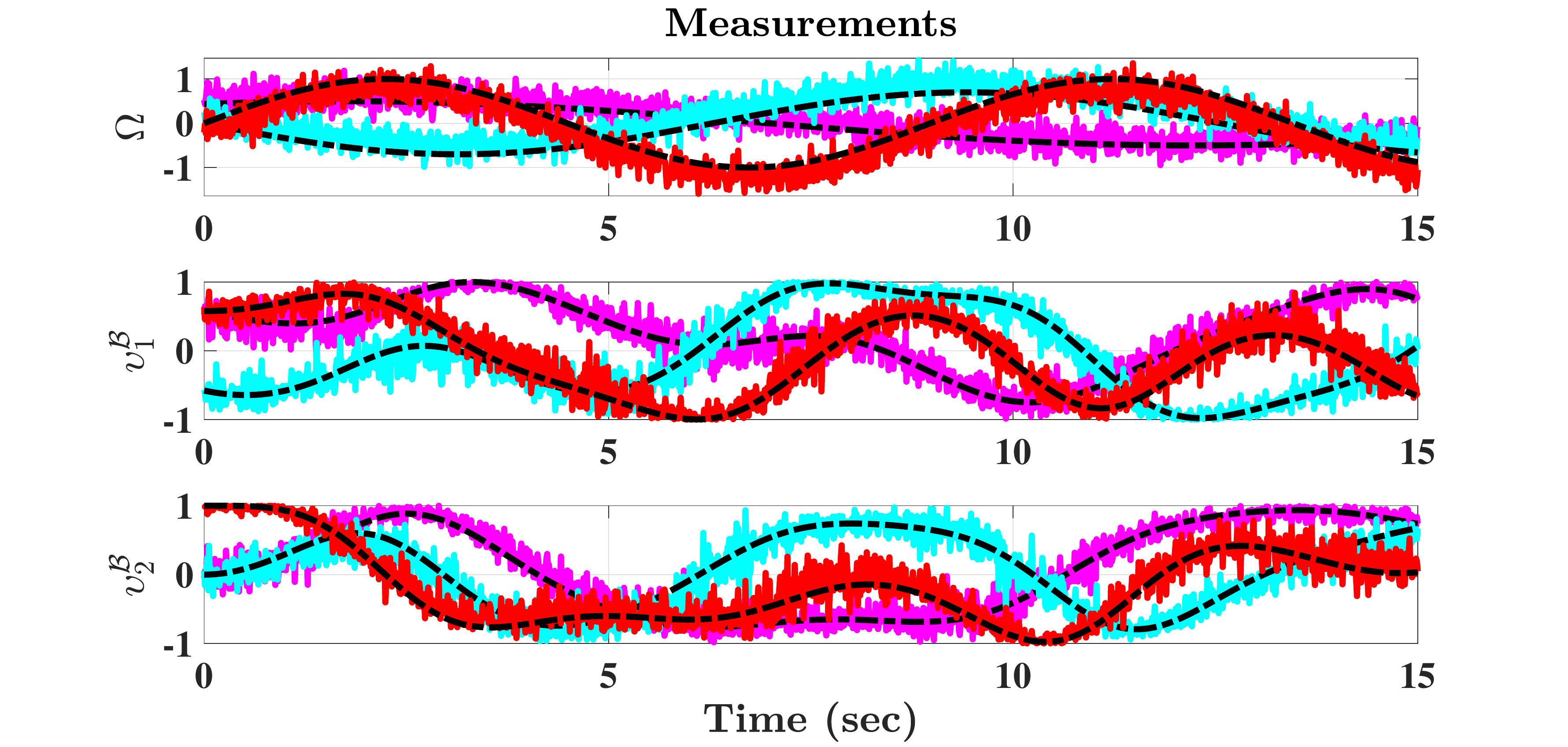}\caption{True values and measurements of $\Omega$, $\upsilon_{1}^{\mathcal{B}}$,
		and $\upsilon_{2}^{\mathcal{B}}$.}
	\label{fig:SO3STCH_1} 
\end{figure}

\begin{figure}[h!]
	\centering{}\includegraphics[scale=0.26]{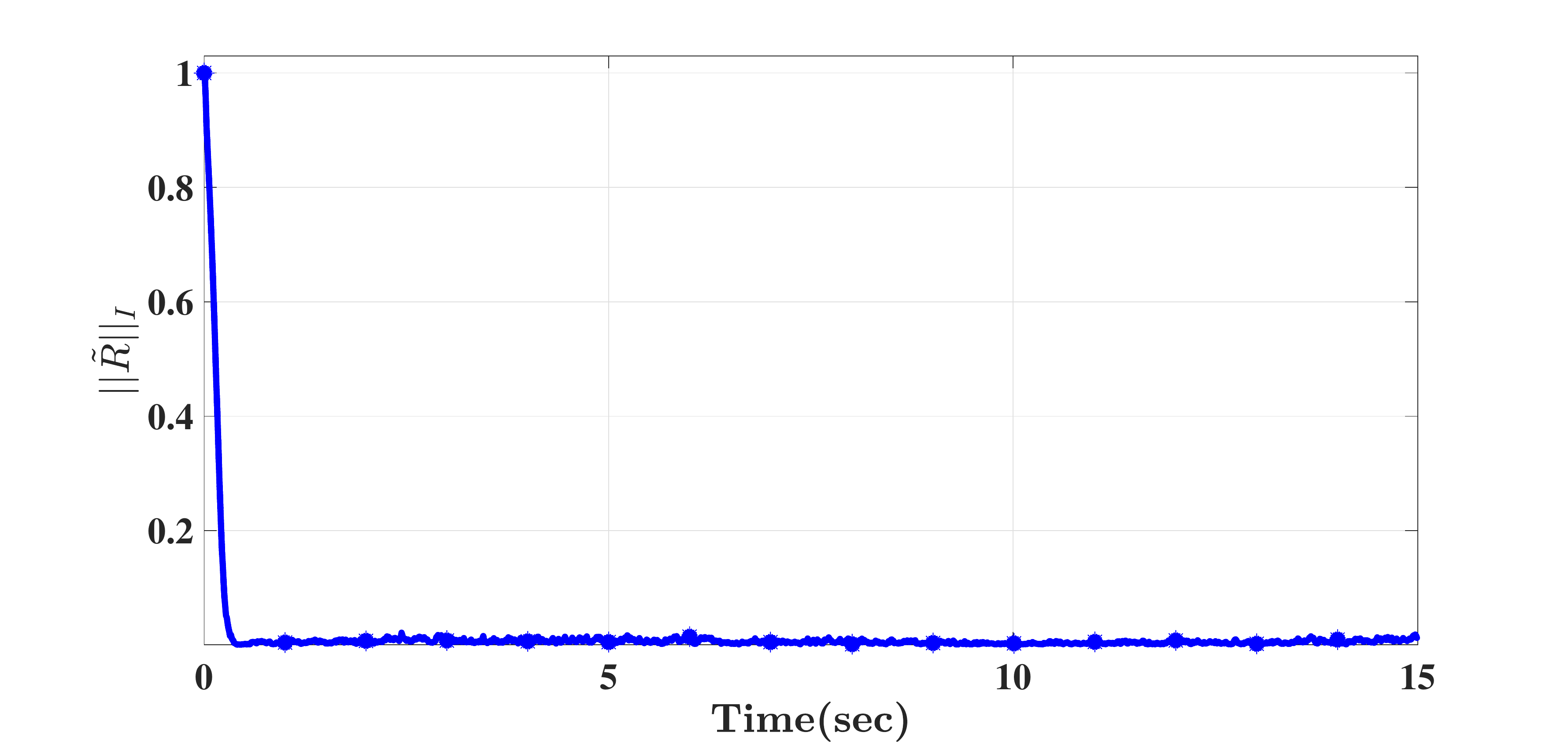}\caption{Tracking performance of normalized Euclidean distance error.}
	\label{fig:SO3STCH_3} 
\end{figure}

\begin{figure*}[h!]
	\centering{}\includegraphics[scale=0.42]{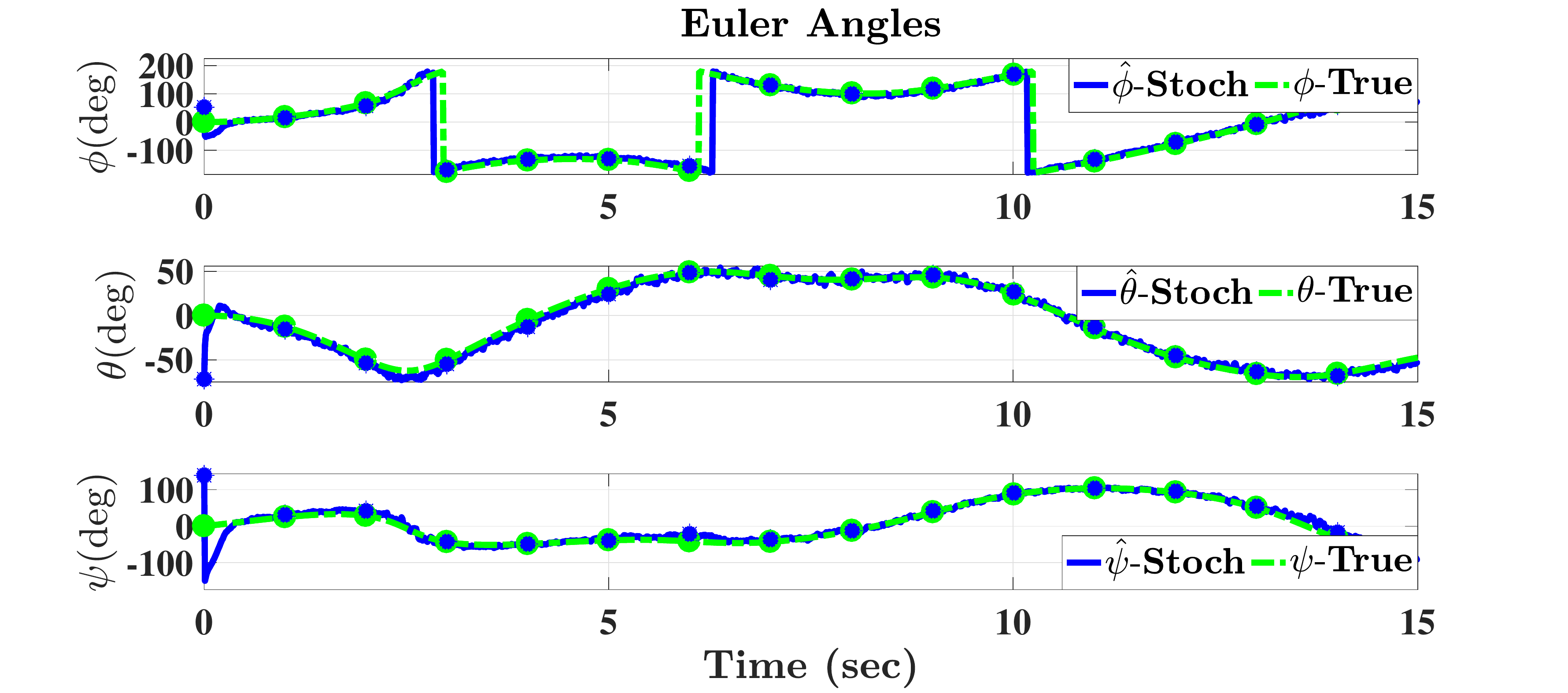}\caption{Tracking performance of Euler angles, proposed filter performance
		vs true trajectories.}
	\label{fig:SO3STCH_5} 
\end{figure*}

Fig. \ref{fig:SO3STCH_1} presents the true angular velocity $\left(\Omega\right)$
and true body-frame vectors as black centerlines and the associated
high values of noise and bias components are represented by colored
solid lines. The robustness of the filter against large initialization
error and high values of noise and bias components is demonstrated
in Fig. \ref{fig:SO3STCH_3}. The normalized Euclidean distance error
$||\tilde{R}||_{I}$ was initiated very close to the unstable equilibria
$\left(+1\right)$, eventually reduced to the neighborhood of the
origin in probability as illustrated in Fig. \ref{fig:SO3STCH_3}. Fig. \ref{fig:SO3STCH_5} illustrate the tracking performance of Euler angles, proposed filter performance
vs true trajectories.

\section{Conclusion \label{sec:SO3STCH_EXP_Conclusion}}

An explicit stochastic nonlinear attitude filter is proposed on $\mathbb{SO}\left(3\right)$.
The proposed filter shares its structure with previously developed
deterministic filters, but in stochastic sense. An alternate attitude
potential function which has not been considered in literature, is
introduced in this work. The resulting stochastic filter ensures that
the errors in Rodriguez vector and estimates are semi-globally uniformly
ultimately bounded in mean square. Numerical results show high convergence
capabilities when large error is initialized in the attitude and high
levels of noise and bias are observed in the vector measurements.

\section*{Appendix A \label{sec:SO3STCH_EXPL_AppendixA} }
\begin{center}
	\textbf{{}Proof of Lemma \ref{Lem:SO3STCH_EXPL_1}} 
	\par\end{center}

Let the attitude be represented by $R\in\mathbb{SO}\left(3\right)$.
From Section \ref{sec:SO3STCH_EXP_Stochastic-Complementary-Filters}
$\sum_{i=1}^{n}s_{i}=3$ which implies that ${\rm Tr}\left\{ M^{\mathcal{I}}\right\} =3$
and the normalized Euclidean distance of $M^{\mathcal{I}}R$ is $\left\Vert M^{\mathcal{I}}R\right\Vert _{I}=\frac{1}{4}{\rm Tr}\left\{ M^{\mathcal{I}}\left(\mathbf{I}_{3}-R\right)\right\} $.
According to angle-axis parameterization in \eqref{eq:SO3STCH_EXPL_att_ang},
one obtains
\begin{align}
\left\Vert M^{\mathcal{I}}R\right\Vert _{I} & =\frac{1}{4}{\rm Tr}\left\{ -M^{\mathcal{I}}\left(\sin(\theta)\left[u\right]_{\times}+\left(1-\cos(\theta)\right)\left[u\right]_{\times}^{2}\right)\right\} \nonumber \\
& =-\frac{1}{4}{\rm Tr}\left\{ M^{\mathcal{I}}\left(1-\cos(\theta)\right)\left[u\right]_{\times}^{2}\right\} \label{eq:SO3STCH_EXPL_append3}
\end{align}
where ${\rm Tr}\left\{ M^{\mathcal{I}}\left[u\right]_{\times}\right\} =0$
as given in identity \eqref{eq:SO3STCH_EXPL_Identity5}. One has \cite{murray1994mathematical}
\begin{equation}
\left\Vert R\right\Vert _{I}=\frac{1}{4}{\rm Tr}\left\{ \mathbf{I}_{3}-R\right\} ={\rm sin}^{2}\left(\theta/2\right)\label{eq:SO3STCH_EXPL_append4}
\end{equation}
and the Rodriguez parameters vector in terms of angle-axis parameterization
is \cite{shuster1981three} $u={\rm cot}\left(\theta/2\right)\rho$.
From identity \eqref{eq:SO3STCH_EXPL_Identity3} $\left[u\right]_{\times}^{2}=-\left\Vert u\right\Vert ^{2}\mathbf{I}_{3}+uu^{\top}$,
the expression in \eqref{eq:SO3STCH_EXPL_append3} becomes
\begin{align*}
\left\Vert M^{\mathcal{I}}R\right\Vert _{I} & =\frac{1}{2}\left\Vert R\right\Vert _{I}u^{\top}\bar{\mathbf{M}}^{\mathcal{I}}u=\frac{1}{2}\left\Vert R\right\Vert _{I}{\rm cot}^{2}\left(\frac{\theta}{2}\right)\rho^{\top}\bar{\mathbf{M}}^{\mathcal{I}}\rho
\end{align*}
From \eqref{eq:SO3STCH_EXPL_append4}, one can find ${\rm cos}^{2}\left(\frac{\theta}{2}\right)=1-\left\Vert R\right\Vert _{I}$
which means
\[
{\rm tan}^{2}\left(\frac{\theta}{2}\right)=\frac{\left\Vert R\right\Vert _{I}}{1-\left\Vert R\right\Vert _{I}}
\]
Consequently, the normalized Euclidean distance is defined in the
sense of Rodriguez parameters vector as
\begin{align}
\left\Vert M^{\mathcal{I}}R\right\Vert _{I} & =\frac{1}{2}\left(1-\left\Vert R\right\Vert _{I}\right)\rho^{\top}\bar{\mathbf{M}}^{\mathcal{I}}\rho=\frac{1}{2}\frac{\rho^{\top}\bar{\mathbf{M}}^{\mathcal{I}}\rho}{1+\left\Vert \rho\right\Vert ^{2}}\label{eq:SO3STCH_EXPL_append_MBR_I}
\end{align}
This proves \eqref{eq:SO3STCH_EXPL_lemm1_1}. The anti-symmetric projection
operator can be defined in terms of Rodriquez parameters vector with
the aid of identity \eqref{eq:SO3STCH_EXPL_Identity1} and \eqref{eq:SO3STCH_EXPL_Identity4}
by
\begin{align*}
\boldsymbol{\mathcal{P}}_{a}\left(M^{\mathcal{I}}R\right)= & \frac{M^{\mathcal{I}}\rho\rho^{\top}-\rho\rho^{\top}M^{\mathcal{I}}+M^{\mathcal{I}}\left[\rho\right]_{\times}+\left[\rho\right]_{\times}M^{\mathcal{I}}}{1+\left\Vert \rho\right\Vert ^{2}}\\
= & \frac{\left[\left({\rm Tr}\left\{ M^{\mathcal{I}}\right\} \mathbf{I}_{3}-M^{\mathcal{I}}+\left[\rho\right]_{\times}M^{\mathcal{I}}\right)\rho\right]_{\times}}{1+\left\Vert \rho\right\Vert ^{2}}
\end{align*}
It follows that the vex operator of the above expression is{\small{}
	\begin{align}
	\Phi\left(M^{\mathcal{I}}R\right)=\mathcal{\mathbf{vex}}\left(\boldsymbol{\mathcal{P}}_{a}\left(M^{\mathcal{I}}R\right)\right) & =\frac{\left(\mathbf{I}_{3}+\left[\rho\right]_{\times}\right)^{\top}}{1+\left\Vert \rho\right\Vert ^{2}}\bar{\mathbf{M}}^{\mathcal{I}}\rho\label{eq:SO3STCH_EXPL_append_MBR_VEX}
	\end{align}
}This shows \eqref{eq:SO3STCH_EXPL_lemm1_2}. The 2-norm of \eqref{eq:SO3STCH_EXPL_append_MBR_VEX}
can be obtained by
\begin{align*}
\left\Vert \Phi\left(M^{\mathcal{I}}R\right)\right\Vert ^{2} & =\frac{\rho^{\top}\bar{\mathbf{M}}^{\mathcal{I}}\left(\mathbf{I}_{3}-\left[\rho\right]_{\times}^{2}\right)\bar{\mathbf{M}}^{\mathcal{I}}\rho}{\left(1+\left\Vert \rho\right\Vert ^{2}\right)^{2}}
\end{align*}
with the aid of identity \eqref{eq:SO3STCH_EXPL_Identity3}, one obtains
\begin{align}
\left\Vert \Phi\left(M^{\mathcal{I}}R\right)\right\Vert ^{2} & =\frac{\rho^{\top}\bar{\mathbf{M}}^{\mathcal{I}}\left(\mathbf{I}_{3}-\left[\rho\right]_{\times}^{2}\right)\bar{\mathbf{M}}^{\mathcal{I}}\rho}{\left(1+\left\Vert \rho\right\Vert ^{2}\right)^{2}}\nonumber \\
& =\frac{\rho^{\top}\left(\bar{\mathbf{M}}^{\mathcal{I}}\right)^{2}\rho}{1+\left\Vert \rho\right\Vert ^{2}}-\frac{\left(\rho^{\top}\bar{\mathbf{M}}^{\mathcal{I}}\rho\right)^{2}}{\left(1+\left\Vert \rho\right\Vert ^{2}\right)^{2}}\nonumber \\
& \geq\underline{\lambda}\left(1-\frac{\left\Vert \rho\right\Vert ^{2}}{1+\left\Vert \rho\right\Vert ^{2}}\right)\frac{\rho^{\top}\bar{\mathbf{M}}^{\mathcal{I}}\rho}{1+||\rho||^{2}}\nonumber \\
& \geq2\underline{\lambda}\left(1-\left\Vert R\right\Vert _{I}\right)\left\Vert M^{\mathcal{I}}R\right\Vert _{I}\label{eq:SO3STCH_EXPL_append_VEX_MI2}
\end{align}
where $\underline{\lambda}=\underline{\lambda}\left(\bar{\mathbf{M}}^{\mathcal{I}}\right)$
and $\left\Vert R\right\Vert _{I}=\left\Vert \rho\right\Vert ^{2}/\left(1+\left\Vert \rho\right\Vert ^{2}\right)$
as defined in \eqref{eq:SO3STCH_EXPL_TR2}. It can be found that
\begin{align}
1-\left\Vert R\right\Vert _{I} & =\frac{1}{4}\left(1+{\rm Tr}\left\{ \left(M^{\mathcal{I}}\right)^{-1}M^{\mathcal{I}}R\right\} \right)\label{eq:SO3STCH_EXPL_append_rho2}
\end{align}
Therefore, from \eqref{eq:SO3STCH_EXPL_append_VEX_MI2}, and \eqref{eq:SO3STCH_EXPL_append_rho2}
the following inequality holds
\begin{align*}
||\Phi\left(M^{\mathcal{I}}R\right)||^{2} & \geq\frac{\underline{\lambda}}{2}\left(1+{\rm Tr}\left\{ \left(M^{\mathcal{I}}\right)^{-1}M^{\mathcal{I}}R\right\} \right)\left\Vert M^{\mathcal{I}}R\right\Vert _{I}
\end{align*}
which proves \eqref{eq:SO3STCH_EXPL_lemm1_4} in Lemma \ref{Lem:SO3STCH_EXPL_1}.

\section*{Acknowledgment}

The authors would like to thank University of Western Ontario for
the funding that made this research possible. Also, the authors would
like to thank \textbf{Maria Shaposhnikova} for proofreading the article.

\bibliographystyle{IEEEtran}
\bibliography{bib_STCH_EXP_SO3}

\begin{thebibliography}{10}
\providecommand{\url}[1]{#1}
\csname url@samestyle\endcsname
\providecommand{\newblock}{\relax}
\providecommand{\bibinfo}[2]{#2}
\providecommand{\BIBentrySTDinterwordspacing}{\spaceskip=0pt\relax}
\providecommand{\BIBentryALTinterwordstretchfactor}{4}
\providecommand{\BIBentryALTinterwordspacing}{\spaceskip=\fontdimen2\font plus
\BIBentryALTinterwordstretchfactor\fontdimen3\font minus
  \fontdimen4\font\relax}
\providecommand{\BIBforeignlanguage}[2]{{%
\expandafter\ifx\csname l@#1\endcsname\relax
\typeout{** WARNING: IEEEtran.bst: No hyphenation pattern has been}%
\typeout{** loaded for the language `#1'. Using the pattern for}%
\typeout{** the default language instead.}%
\else
\language=\csname l@#1\endcsname
\fi
#2}}
\providecommand{\BIBdecl}{\relax}
\BIBdecl

\bibitem{mahony2008nonlinear}
R.~Mahony, T.~Hamel, and J.-M. Pflimlin, ``Nonlinear complementary filters on
  the special orthogonal group,'' \emph{IEEE Transactions on Automatic
  Control}, vol.~53, no.~5, pp. 1203--1218, 2008.

\bibitem{crassidis2007survey}
J.~L. Crassidis, F.~L. Markley, and Y.~Cheng, ``Survey of nonlinear attitude
  estimation methods,'' \emph{Journal of guidance, control, and dynamics},
  vol.~30, no.~1, pp. 12--28, 2007.

\bibitem{choukroun2006novel}
D.~Choukroun, I.~Y. Bar-Itzhack, and Y.~Oshman, ``Novel quaternion kalman
  filter,'' \emph{IEEE Transactions on Aerospace and Electronic Systems},
  vol.~42, no.~1, pp. 174--190, 2006.

\bibitem{lefferts1982kalman}
E.~J. Lefferts, F.~L. Markley, and M.~D. Shuster, ``Kalman filtering for
  spacecraft attitude estimation,'' \emph{Journal of Guidance, Control, and
  Dynamics}, vol.~5, no.~5, pp. 417--429, 1982.

\bibitem{markley2003attitude}
F.~L. Markley, ``Attitude error representations for kalman filtering,''
  \emph{Journal of guidance, control, and dynamics}, vol.~26, no.~2, pp.
  311--317, 2003.

\bibitem{crassidis2003unscentedy}
J.~L. Crassidis and F.~L. Markley, ``Unscented filtering for spacecraft
  attitude estimation,'' \emph{Journal of guidance, control, and dynamics},
  vol.~26, no.~4, pp. 536--542, 2003.

\bibitem{Arulampalam2002Particle}
M.~S. Arulampalam, S.~Maskell, N.~Gordon, and T.~Clapp, ``A tutorial on
  particle filters for online nonlinear/non-gaussian bayesian tracking,''
  \emph{IEEE Transactions on Signal Processing}, vol.~50, no.~2, pp. 174--185,
  2002.

\bibitem{shuster1981three}
M.~D. Shuster and S.~D. Oh, ``Three-axis attitude determination from vector
  observations,'' \emph{Journal of Guidance, Control, and Dynamics}, vol.~4,
  pp. 70--77, 1981.

\bibitem{zlotnik2017nonlinear}
D.~E. Zlotnik and J.~R. Forbes, ``Nonlinear estimator design on the special
  orthogonal group using vector measurements directly,'' \emph{IEEE
  Transactions on Automatic Control}, vol.~62, no.~1, pp. 149--160, 2017.

\bibitem{grip2012attitude}
H.~F. Grip, T.~I. Fossen, T.~A. Johansen, and A.~Saberi, ``Attitude estimation
  using biased gyro and vector measurements with time-varying reference
  vectors,'' \emph{IEEE Transactions on Automatic Control}, vol.~57, no.~5, pp.
  1332--1338, 2012.

\bibitem{hashim2017adaptive}
H.~A. Hashim, S.~El-Ferik, and F.~L. Lewis, ``Adaptive synchronisation of
  unknown nonlinear networked systems with prescribed performance,''
  \emph{International Journal of Systems Science}, vol.~48, no.~4, pp.
  885--898, 2017.

\bibitem{hashim2017neuro}
H.~A. Hashim, S.~El-Ferik, and F.~L. Lewis, ``Neuro-adaptive cooperative
  tracking control with prescribed performance of unknown higher-order
  nonlinear multi-agent systems,'' \emph{International Journal of Control}, pp.
  1--16, 2017.

\bibitem{shuster1993survey}
M.~D. Shuster, ``A survey of attitude representations,'' \emph{Navigation},
  vol.~8, no.~9, pp. 439--517, 1993.

\bibitem{khasminskii1980stochastic}
R.~Khasminskii, \emph{Stochastic stability of differential equations}.\hskip
  1em plus 0.5em minus 0.4em\relax Rockville, MD: S \& N International, 1980.

\bibitem{hashim2018SO3Stochastic}
H.~A. Hashim, L.~J. Brown, and K.~McIsaac, ``Nonlinear stochastic attitude
  filters on the special orthogonal group 3: Ito and stratonovich,'' \emph{IEEE
  Transactions on Systems, Man, and Cybernetics: Systems}, pp. 1--13, 2018
  (Submitted).

\bibitem{ito1984lectures}
K.~Ito and K.~M. Rao, \emph{Lectures on stochastic processes}.\hskip 1em plus
  0.5em minus 0.4em\relax Tata institute of fundamental research, 1984,
  vol.~24.

\bibitem{deng2001stabilization}
H.~Deng, M.~Krstic, and R.~J. Williams, ``Stabilization of stochastic nonlinear
  systems driven by noise of unknown covariance,'' \emph{IEEE Transactions on
  Automatic Control}, vol.~46, no.~8, pp. 1237--1253, 2001.

\bibitem{ji2006adaptive}
H.-B. Ji and H.-S. Xi, ``Adaptive output-feedback tracking of stochastic
  nonlinear systems,'' \emph{IEEE Transactions on Automatic Control}, vol.~51,
  no.~2, pp. 355--360, 2006.

\bibitem{hashim2018SE3Stochastic}
H.~A. Hashim, L.~J. Brown, and K.~McIsaac, ``Nonlinear stochastic position and
  attitude filter on the special euclidean group 3,'' \emph{Journal of the
  Franklin Institute}, pp. 1--27, 2018 (Submitted).

\bibitem{bullo2004geometric}
F.~Bullo and A.~D. Lewis, \emph{Geometric control of mechanical systems:
  modeling, analysis, and design for simple mechanical control systems}.\hskip
  1em plus 0.5em minus 0.4em\relax Springer Science \& Business Media, 2004,
  vol.~49.

\bibitem{murray1994mathematical}
R.~M. Murray, Z.~Li, S.~S. Sastry, and S.~S. Sastry, \emph{A mathematical
  introduction to robotic manipulation}.\hskip 1em plus 0.5em minus 0.4em\relax
  CRC press, 1994.

\end{thebibliography}

\end{document}